\documentclass[12pt]{amsart}
\usepackage[T1]{fontenc}
\usepackage[utf8]{inputenc}
\usepackage{microtype}
\usepackage{amsmath,amssymb,amsthm,mathrsfs}
\usepackage{dsfont}
\usepackage{esint}
\usepackage{slashed}

\usepackage{mathtools}

\mathtoolsset{showonlyrefs}

\usepackage{thmtools, thm-restate}

\usepackage{tikz-cd}

\tikzset{
	pf/.style={commutative diagrams/.cd, every arrow, every label},
	surj/.style=commutative diagrams/two heads,
	inj/.style=commutative diagrams/hook,
	gl/.style=commutative diagrams/equal,
	mat/.style={matrix of math nodes, commutative diagrams/.cd, every cell},
	dr/.style={matrix of math nodes, commutative diagrams/.cd, every cell, column sep=small},
	seq/.style={matrix of math nodes, commutative diagrams/.cd, every cell,column sep=small}
}

\newenvironment{diag*}
	{\[\begin{tikzpicture}[commutative diagrams/.cd, every diagram]}
	{\end{tikzpicture}\]\ignorespacesafterend}

\newenvironment{diag}
	{\begin{equation}\begin{tikzpicture}[commutative diagrams/.cd, every diagram]}
	{\end{tikzpicture}\end{equation}\ignorespacesafterend}

\usepackage{hyperref}

\usepackage[margin=0.9in]{geometry}

\newcommand{\al}{\alpha}
\newcommand{\be}{\beta}
\newcommand{\p}{\partial}
\newcommand{\R}{\mathbb{R}}
\newcommand{\na}{\nabla}

\newcommand{\la}{\langle}
\newcommand{\ra}{\rangle}

\newcommand{\D}{\slashed{D}}
\newcommand{\pd}{\slashed{\partial}}
\newcommand{\ACM}{J_M}
\newcommand{\ACSBP}{J_\SpinorBundleP}
\newcommand{\ACSBMI}{I_\SpinorBundleM}
\newcommand{\ACSBMJ}{J_\SpinorBundleM}
\newcommand{\ACSBMK}{K_\SpinorBundleM}
\newcommand{\quaternions}{\mathbb{H}}
\newcommand{\CliffordBundleM}{{\operatorname{Cl}\left(M, -g\right)}}
\newcommand{\CliffordBundleP}{{\operatorname{Cl}\left(M, g\right)}}
\newcommand{\ComplexHalfSpinorBundle}{W}
\newcommand{\DiracOperatorM}{\pd}
\newcommand{\DiracOperatorP}{{\pd_{\SpinorBundleP}}}
\newcommand{\DiracOperatorPC}{{\pd_{\SpinorBundleP}^\mathbb{C}}}
\newcommand{\CliffordAlgebraM}{\operatorname{Cl}_{0,2}}
\newcommand{\CliffordAlgebraP}{\operatorname{Cl}_{2,0}}
\newcommand{\SpinGroup}{\operatorname{Spin}(2)}
\newcommand{\SpinorBundleM}{S}
\newcommand{\SpinorBundleMEven}{\SpinorBundleM^0}
\newcommand{\SpinorBundleMOdd}{\SpinorBundleM^1}
\newcommand{\SpinorBundleP}{\Sigma}
\newcommand{\SpinorBundlePC}{\Sigma^{\mathbb{C}}}
\newcommand{\SpinorHermitianFormPC}{h_{\SpinorBundlePC}}
\newcommand{\SpinorHermitianFormM}{h_{\SpinorBundleM}}
\newcommand{\SpinorMetricM}{g_{\SpinorBundleM}}
\newcommand{\SpinorMetricP}{g_{\SpinorBundleP}}
\newcommand{\SpinStructure}{P_{\SpinGroup}}
\newcommand{\vph}{\varphi}
\newcommand{\A}{\mathbb{A}}

\newcommand{\MS}[1]{M^{#1}}

\newcommand{\dd}{\mathop{}\!\mathrm{d}}

\DeclareMathOperator{\Rm}{R}

\DeclareMathOperator{\diverg}{div}

\DeclareMathOperator{\Id}{Id}

\newtheorem{thm}{Theorem}[section]
\newtheorem{Def}[thm]{Definition}
\newtheorem{lemma}[thm]{Lemma}

\newtheorem{rmk}[thm]{Remark}

\title{Regularity of Solutions of the Nonlinear Sigma Model with Gravitino}

\begin{document}

\author{Jürgen Jost, Enno Keßler, Jürgen Tolksdorf, Ruijun Wu, Miaomiao Zhu}

\address{Max Planck Institute for Mathematics in the Sciences\\Inselstr. 22--26\\D-04103 Leipzig, Germany}
	\email{jjost@mis.mpg.de}

\address{Max Planck Institute for Mathematics in the Sciences\\Inselstr. 22--26\\D-04103 Leipzig, Germany}
	\email{Juergen.Tolksdorf@mis.mpg.de}

\address{Max Planck Institute for Mathematics in the Sciences\\Inselstr. 22--26\\D-04103 Leipzig, Germany}
	\email{Enno.Kessler@mis.mpg.de}

\address{Max Planck Institute for Mathematics in the Sciences\\Inselstr. 22--26\\D-04103 Leipzig, Germany}
	\email{Ruijun.Wu@mis.mpg.de}

\address{School of Mathematical Sciences, Shanghai Jiao Tong University\\Dongchuan Road 800\\200240 Shanghai, P.R.China}
	\email{mizhu@sjtu.edu.cn}

\thanks{%
\emph{Acknowledgements:}
Ruijun Wu thanks the International Max Planck Research School Mathematics in the Sciences for financial support.
Miaomiao Zhu was supported in part by the National Natural Science Foundation of China (No. 11601325).
}

\date{\today}

\begin{abstract}
	We propose a geometric setup to study analytic aspects of a variant of the super symmetric two-dimensional nonlinear sigma model.
	This functional extends the functional of Dirac-harmonic maps by gravitino fields.
	The system of Euler--Lagrange equations of the two-dimensional nonlinear sigma model with gravitino is calculated explicitly.
	The gravitino terms pose additional analytic difficulties to show smoothness of its weak solutions which are overcome using Rivière's regularity theory and Riesz potential theory.
\end{abstract}


\maketitle

\section{Introduction}
The various versions of the two-dimensional sigma models are among the most important and best studied models of quantum field theory.
On one hand, such models possess important symmetries, in particular conformal invariance.
On the other hand, they can be analyzed in detail with difficult, but currently available mathematical methods.
Here, we shall investigate its probably most general and physically and mathematically richest version, the two-dimensional supersymmetric nonlinear sigma-model, introduced in~\cite{brink1976locally, deser1976complete}.
This model possesses a subtle mathematical structure, see~\cite{deligne1999quantum, jost2009geometry}.
The physical and mathematical structure of the model depends on the symmetries it possesses.
These include generalized conformal invariance, super Weyl symmetry, and supersymmetry, hence the name of the model.

While supersymmetry requires anti-commuting variables, a version of this model with all fields commuting has been intensively studied by mathematicians in the last decade.
The mathematical analysis started with various reduced forms of this model.
The simplest instance are harmonic functions, which correspond to the linear sigma model,
and they have played an important role in analysis and geometry for a long time.
The nonlinear version leads to harmonic maps instead of functions,
and these are likewise well studied objects with many applications in geometric analysis.
In the super version, the map gets coupled with a super partner, a vector spinor.
Chen--Jost--Li--Wang~\cite{chen2006dirac, chen2005regularity} initiated the analysis of such coupled fields, which they called Dirac-harmonic maps.
The full physical model contains still more additional terms, some of which were considered in~\cite{chen2007, branding2015some, branding2015energy, branding2016, jost2015geometric}.
Based on those works, we are now in a position to address the full model, including the gravitino terms.
The supersymmetric action functional has been mathematically studied from an algebraic and geometric perspective in a systematic way in~\cite{jost2014super}.
Here we shall start to explore the analytic aspects.

Let $(M,g)$ be a closed, oriented surface and $(N,h)$ a closed Riemannian manifold.
We will study the super action functional $\A$ defined on the space
\begin{equation}
	\mathcal{X}^{1,2}_{1,4/3}(M,N)=\{(\phi,\psi)\big| \phi\in W^{1,2}(M,N), \psi\in\Gamma^{1,4/3}(S\otimes\phi^*TN)\},
\end{equation}
where by $\Gamma^{1,4/3}(S\otimes\phi^*TN)$ we mean the space of $W^{1,4/3}$ sections of the twisted spinor bundle \(S\otimes \phi^*TN\).
Furthermore, in this paper the Riemannian metric \(g\) and the gravitino \(\chi\) are considered parameters of the functional.
Even though an \(L^4\)-integrability condition suffices for the finiteness of \(\A\), we will always assume the gravitino \(\chi\) is a smooth section of \(\SpinorBundleM\otimes TM\).
The action functional is
\begin{equation}%
\label{A-intro}
	\begin{split}
		\A(\phi, \psi;g, \chi)\coloneqq \int_M & |\dd \phi|_{T^*M\otimes \phi^*TN}^2
			+ \langle \psi, \D \psi \rangle_{S\otimes \phi^*TN}  \\
		&  -4\langle (\mathds{1}\otimes\phi_*)(Q\chi), \psi \rangle_{S\otimes\phi^*TN}
			-|Q\chi|^2_{S\otimes TM} |\psi|^2_{S\otimes \phi^*TN}
			-\frac{1}{6} \Rm^{N}(\psi) \dd vol_g.
	\end{split}
\end{equation}
Here $Q$ is a projection operator mapping to a subspace of \(\SpinorBundleM\otimes TM\), \(\mathds{1}\otimes\phi_*\colon \SpinorBundleM\otimes TM\to \SpinorBundleM\otimes \phi^*TN\) and  \(R^N(\psi)\) is a contraction of the pullback of the curvature of $N$ along $\phi$ with the field \(\psi\) to the fourth order.

While the precise geometric setup will be explained in Section~\ref{Sec:Preliminaries}, we already give local expressions for the third and fifth summand.
Let \(\{e_\alpha\}\) be a local orthonormal frame of \(TM\) and \(\{y^i\}\) local coordinates on \(N\).
Writing \(\chi = \chi^\alpha\otimes e_\alpha\) and \(\psi = \psi^i\otimes \phi^*\left(\partial_{y^i}\right)\) it holds
\begin{align}
	-4\langle (\mathds{1}\otimes\phi_*)(Q\chi), \psi \rangle_{S\otimes\phi^*TN}
		&=2\langle e_\al \cdot e_\be \cdot \chi^\al \otimes \phi_* e_\be, \psi \rangle_{S\otimes \phi^*TN}, \\
	-\frac{1}{6}\Rm^{N}(\psi)
		&=-\frac{1}{6}\Rm^{N}_{ijkl}\langle \psi^i, \psi^k \rangle_S \langle \psi^j,\psi^l \rangle_S.
\end{align}

Since the action functional is somewhat involved and contains many different fields and at the same time possesses rich symmetries, the derivation of the associated Euler--Lagrange equations requires substantial computations.
This will be the first achievement of this paper.
The result is:
\begin{thm}[restate=ELThm, label=thm:EL]
	The Euler--Lagrange equations for the super action functional $\A$ are given by
	\begin{equation} \label{EL-eq}
		\begin{split}
			\tau(\phi)=&\frac{1}{2}\Rm^{\phi^*TN}(\psi, e_\al\cdot\psi)\phi_* e_\al-\frac{1}{12}S\na R(\psi)   \\
			&  -(\langle \na^S_{e_\be}(e_\al \cdot e_\be \cdot \chi^\al), \psi \rangle_S
		+ \langle e_\al \cdot e_\be \cdot \chi^\al, \na^{S\otimes\phi^*TN}_{e_\be} \psi \rangle_S),  \\
		\D\psi =& |Q\chi|^2\psi +\frac{1}{3}SR(\psi)+2(\mathds{1}\otimes \phi_*)Q\chi.
	\end{split}
\end{equation}
\end{thm}
These equations already make the growth order transparent with which the various fields enter.
$SR(\psi)$ stands for a term involving the curvature of the target $N$ that is cubic in $\psi$, see~\eqref{def-SR} and $S\na R(\psi)$ involves derivatives of that curvature and is quartic in $\psi$, see~\eqref{def-SNR}.

We shall then turn to the properties of their solutions.
More precisely, we want to show the regularity of weak solutions, that is, those that satisfy the Euler--Lagrange equations in the sense of distributions.

The basic issues in geometric analysis are the existence, uniqueness and smoothness of nontrivial critical points.
That is, one wishes to show the existence of weak solutions and then their uniqueness and regularity.
In this paper, we settle the smoothness.
The Euler--Lagrange equations~\eqref{EL-eq} of this action functional turn out to be critical for the Sobolev framework, in the sense that, with initial data assumed to lie in some Sobolev spaces, the classical bootstrap arguments are not strong enough to improve the regularity.
That is, the powerful scheme of elliptic regularity theory does not directly apply,
and we need to utilize the structure of the equations, and in particular their symmetries, in a subtler way.
Our analytical tools are the Morrey spaces, which can be viewed as finer subspaces of the Lebesgue spaces.
With estimates on Riesz potentials, we can then iteratively improve the regularity, and get the system away from the critical case.
Related methods have been used in~\cite{wang2010remark, sharp2016regularity, branding2015some}.
Then the Rivi\`ere regularity theory (see e.g.~\cite{riviere, riviere2010conformally, sharp2013decay}) can be applied to the map component of the critical pairs.
Finally, we can show that
\begin{thm}%
\label{theorem 1}
	The critical points of the super action functional
	\begin{equation}
		\begin{split}
			\A\colon\mathcal{X}^{1,2}_{1,4/3}(M,N) & \to \R, \\
			(\phi,\psi)& \mapsto \A(\phi,\psi;g,\chi),
		\end{split}
	\end{equation}
	are smooth, provided \(g\) and \(\chi\) are smooth.
\end{thm}
This result should also help in finding solutions of its associated Euler--Lagrange equations.
Moreover, our method is of interest in its own right, as we shall explain later.
Further geometric and analytic aspects of this model will be addressed in subsequent work.

As in the aforementioned works, we shall work with the version of the model that only has commuting fields.
As explained in~\cite{chen2011boundary}, this depends on an appropriate representation of the Clifford algebra involved.
Thus, in contrast to~\cite{jost2014super}, we shall not have to work in the category of supermanifolds,
but can confine ourselves to the setting of Riemannian geometry.
Yet, in the framework of supermanifolds, the action functional~\eqref{A-intro} and its symmetries obtain a natural geometric interpretation.
In~\cite{jost2014super} it was shown that the fields \(g\) and \(\chi\) determine a \emph{super Riemann surface}, a super geometric generalization of a Riemann surface.
Recall that Teichmüller theory can be developed with the help of the harmonic action functional.
The functional \(\A\) can be seen as a super analogue of the harmonic action functional.
Hence it is expected that an understanding of the solution space of the Euler--Lagrange equations~\eqref{EL-eq} helps to study geometric properties of the moduli space of super Riemann surfaces.

Concerning the organization of this paper, we shall first set up the geometric background for the model and introduce the action functional as well as its basic properties.
Then we shall derive its Euler--Lagrange equations.
For our regularity scheme, we need to bring the equations into a suitable form.
This treatment of the Euler--Lagrange equations which builds upon~\cite{zhu2009regularity, wang2009regularity, chen2011boundary, sharp2016regularity, branding2015some} is crucial for our paper, and we hope that it will also be useful for the further mathematical investigation of the model.
We can then finally show the regularity of weak solutions of the Euler--Lagrange equations.
The main lemma in improving the regularity appears in the last section in a somewhat more general form than needed for our present purposes.

\section{Preliminaries}%
\label{Sec:Preliminaries}
In this section we summarize the geometrical background and thereby also fix the notation used in what follows in the subsequent sections.
The main purpose of this section is to provide a geometrical setup such that the action functional~\eqref{A-intro} can be seen as a real-valued action functional with non-vanishing Dirac-action.
Those two requirements will be satisfied using a real four-dimensional spinor representation.
In contrast, in the description of non-linear sigma models on two-dimensional manifolds, two-dimensional real or complex spinor representations are usually taken into account, see for example~\cite{chen2006dirac, jost2014super}.
For the convenience of the reader we add some comments on how these different geometrical settings are related.

\subsection{}
Let \((M,g)\) be a closed, oriented, two-dimensional Riemannian spin manifold with fixed spin structure.
The corresponding \(\SpinGroup\) principal bundle is denoted by \(\SpinStructure\).
For any bilinear form \(b\) on \(TM\) we denote by \(\operatorname{Cl}(M, b)\) the corresponding Clifford algebra bundle, which is isomorphic to the quotient of the tensor algebra by the two-sided ideal generated by
\begin{equation}
	X\otimes Y + Y\otimes X - 2b(X, Y),
\end{equation}
where \(X, Y\in\Gamma(TM)\).
In the following we will only use \(b=\pm g\).

The typical fiber of \(\CliffordBundleP\), denoted by \(\CliffordAlgebraP\), is a simple algebra and isomorphic to \(\mathfrak{gl}(2,\mathbb{R})\).
We denote this isomorphism by \(\gamma^+\colon \CliffordAlgebraP\to\mathfrak{gl}(2,\mathbb{R})\).
Hence, the spinor bundle of \(\CliffordBundleP\) is given by \(\SpinorBundleP=\SpinStructure\times_{\gamma^+} \mathbb{R}^2\) where \(\SpinGroup\subset\mathfrak{gl}(2,\mathbb{R})\) acts by left-multiplication on \(\mathbb{R}^2\).
We denote the Clifford multiplication of a tangent vector \(X\) with \(s\in\Gamma(\SpinorBundleP)\) by \(\gamma^+(X)s\) or simply by \(X\cdot s\) if no confusion arises.
By its construction as an associated bundle to \(\SpinStructure\), the bundle \(\SpinorBundleP\) possesses a natural fiber metric~\(\SpinorMetricP\) such that the Clifford action by tangent vectors is symmetric.
The Levi-Civita connection on~\(TM\) lifts to the spin connection \(\nabla^{\SpinorBundleP}\) on \(\SpinorBundleP\).

The spin Dirac operator is defined with respect to a local \(g\)-orthonormal frame \(e_a\) by \(\DiracOperatorP s = e_a \cdot \nabla^{\SpinorBundleP}_{e_a} s\) for a section \(s\) of \(\SpinorBundleP\).
It is easy to see that \(\DiracOperatorP\) is antisymmetric and hence for any spinor~\(s\) the Dirac action vanishes, that is,
\begin{equation}%
\label{eq:Diracvanish}
	\int_M \SpinorMetricP\left(s, \DiracOperatorP s\right) \dd{vol}_g = 0.
\end{equation}
In order to avoid the vanishing of the Dirac action one may work with anti-commuting spinors, see for example~\cite{jost2014super} and references therein.
Another possibility to obtain a non-vanishing Dirac action is to consider the complexification \(\SpinorBundlePC=\SpinorBundleP\otimes\mathbb{C}\) and the resulting Hermitian form~\(\SpinorHermitianFormPC\).
Then the operator \(i\DiracOperatorPC\), where \(\DiracOperatorPC\) is the complex linear extension of \(\DiracOperatorP\), is symmetric.
Consequently the Dirac action
\begin{equation}
	\int_M \SpinorHermitianFormPC\left( s, i\DiracOperatorPC s\right) \dd{vol}_g, \qquad s\in\Gamma(\SpinorBundlePC)
\end{equation}
does not vanish identically and is real valued.
An equivalent reformulation of this approach was introduced in~\cite{chen2006dirac}.
Notice, however, that the third summand of~\eqref{A-intro} involves a scalar product of two different spinors.
If this scalar product were to be implemented by \(\SpinorHermitianFormPC\), the action functional~\eqref{A-intro} would not be guaranteed to be real-valued.
Whence we replace the two-dimensional complex spinor representation of the approach presented in~\cite{chen2006dirac} by a four-dimensional real one.
This step will be explained next.

\subsection{}
The typical fiber of the Clifford algebra bundle \(\CliffordBundleM\) is the Clifford algebra \(\CliffordAlgebraM\).
As a real associative algebra with unit the Clifford algebra \(\CliffordAlgebraM\) is isomorphic to the quaternions~\(\quaternions\).
Consequently, the left-regular representation of \(\CliffordAlgebraM\) on itself is irreducible.
Hence, we may regard the vector bundle \(\SpinorBundleM=\SpinStructure\times_{\SpinGroup}\CliffordAlgebraM\) as a spinor bundle, where \(\SpinGroup\subset\CliffordAlgebraM\) acts via the left-regular representation of \(\CliffordAlgebraM\).
The spinor bundle \(\SpinorBundleM\) is a four-dimensional real vector bundle.

Notice that \(\CliffordAlgebraM\) is a \(\mathbb{Z}_2\)-graded module over the \(\mathbb{Z}_2\)-graded algebra \(\CliffordAlgebraM\).
As a consequence also the spinor bundle \(\SpinorBundleM=\SpinorBundleMEven\oplus \SpinorBundleMOdd\) is a \(\mathbb{Z}_2\)-graded module over the \(\mathbb{Z}_2\)-graded algebra bundle \(\CliffordBundleM\).
Here, both the even and the odd part of \(\SpinorBundleM\) are isomorphic to \(\SpinorBundleP\) as associated bundles to \(\SpinStructure\).
The Clifford action \(\gamma(X)\) of a tangent vector \(X\) on \(S\) must be of the form
\begin{equation}
\label{eq:RepresentationM}
	\gamma(X) =
	\begin{pmatrix}
		0 & -\gamma^+(X) \\
		\gamma^+(X) & 0 \\
	\end{pmatrix}
\end{equation}
because it is odd with respect to the \(\mathbb{Z}_2\)-grading.
Recall that \(\gamma^+(X)\) denotes the Clifford multiplication of \(X\) on \(\SpinorBundleP\), where \(X\) is considered as an element of \(\CliffordBundleP\).

The induced metric and spin connection on \(\SpinorBundleM\) are denoted, respectively, by \(\SpinorMetricM = \SpinorMetricP\oplus\SpinorMetricP\) and \(\nabla^\SpinorBundleM = \nabla^\SpinorBundleP\oplus\nabla^\SpinorBundleP\).
The action of \(TM\subset\CliffordBundleM\) on \(\SpinorBundleM\) is skew-symmetric with respect to \(\SpinorMetricM\).
Whence the spin Dirac operator \(\DiracOperatorM = e_\alpha\cdot\nabla^S_{e_\alpha}\colon \Gamma(\SpinorBundleM)\to\Gamma(\SpinorBundleM)\) is symmetric with respect to the \(L^2(\SpinorBundleM)\) scalar product
\begin{equation}
	\left<s, t\right>_{L^2(\SpinorBundleM)} =\int_M \SpinorMetricM(s, t) \dd{vol}_g \qquad s,t\in\Gamma(\SpinorBundleM).
\end{equation}
In particular, the Dirac action \(\left<s,\DiracOperatorM s\right>_{L^2(\SpinorBundleM)}\) is non-trivial, as opposed to its \(\CliffordAlgebraP\) counterpart~\eqref{eq:Diracvanish}.
Furthermore, \(\DiracOperatorM\) is essentially self-adjoint, see~\cite[Chapter II, Theorem 5.7]{lawson1989spin}\footnote{Notice that this reference uses a different sign convention and naming scheme for Clifford algebras.}.

\subsection{}
We now explain the different complex structures on the spinor bundles \(\SpinorBundleP\) and \(\SpinorBundleM\).
This will be needed later on and help to clarify the relation to the geometrical setup introduced in~\cite{chen2006dirac}.

Recall that the Riemann surface \(M\) possesses an integrable almost complex structure \(\ACM\) that is defined by
\begin{equation}
	g(\ACM X, Y) = \dd{vol}_g(X, Y)
\end{equation}
for all tangent vectors \(X\) and \(Y\).
Consequently, the tangent bundle \(TM\) is a holomorphic line bundle.

When seen as \(TM\subset\CliffordBundleP\), the almost complex structure \(\ACM\) can be realized as right-multiplication by the volume form \(\omega\).
With respect to a local oriented \(g\)-orthonormal frame \(e_\alpha\) the volume form is given by \(\omega=e_1\cdot e_2\).
Similarly, left-multiplication by \(\omega\) induces an almost complex structure on \(\SpinorBundleP\), which we denote by \(\ACSBP\).
The bundle \(\SpinorBundlePC=\SpinorBundleP\otimes\mathbb{C}\) decomposes in eigen bundles of \(i\ACSBP^\mathbb{C}\), where \(\ACSBP^\mathbb{C}\) denotes the complex linear extension of \(\ACSBP\).
The complex line bundles \(\ComplexHalfSpinorBundle = (\SpinorBundleP, \ACSBP)\) of eigenvalue \(-1\) and \(\overline{\ComplexHalfSpinorBundle} = (\SpinorBundleP, -\ACSBP)\) of eigenvalue \(+1\) are, respectively, the so-called bundles of ``left- and right-handed'' Weyl spinors.

On \(\ComplexHalfSpinorBundle=(\SpinorBundleP, \ACSBP)\) there is a bilinear form with values in \(T^*M\) given by
\begin{equation}
	\SpinorMetricP(s, e_\alpha\cdot t) e^\alpha, \qquad s,t \in \Gamma(\SpinorBundleP),
\end{equation}
where \(e^\alpha\) is the dual basis to the \(g\)-orthonormal frame \(e_\alpha\).
The compatibility of Clifford multiplication and almost complex structures, \(\left(\ACM X\right)\cdot t = X\cdot \ACSBP t = -\ACSBP\left(X\cdot t\right)\), turns the bilinear form into a complex linear isomorphism \(\ComplexHalfSpinorBundle\otimes_\mathbb{C}\ComplexHalfSpinorBundle=T^*M\).
In particular \(\ComplexHalfSpinorBundle\) is a holomorphic vector bundle.
In other words, holomorphic tangent vector fields on a Riemann surface with fixed spin structure have a ``square root''.
Conversely, on a Riemann surface $(M,\ACM)$ every square root of~\(TM\) gives rise to a spin structure on $M$.

Obviously, the complex vector bundle \((S, \ACSBP\oplus\ACSBP)\) is isomorphic to \(\ComplexHalfSpinorBundle\oplus\ComplexHalfSpinorBundle\).
In addition, the spinor bundle \(\SpinorBundleM\) possesses three almost complex structures \(\ACSBMI, \ACSBMJ, \ACSBMK\in \operatorname{End}(\SpinorBundleM)\) that commute with the Clifford multiplication and satisfy the quaternionic relations: \(\ACSBMI^2 = \ACSBMJ^2 = \ACSBMK^2 = -\Id_\SpinorBundleM\) and \(\ACSBMI = \ACSBMJ\circ \ACSBMK = - \ACSBMK\circ \ACSBMJ\), etc.
Explicitly, they are given by \(\ACSBMI(s,t) = (-t,s)\), \(\ACSBMJ(s,t) = (\ACSBP s,-\ACSBP t)\) and \(\ACSBMK(s,t) = (\ACSBP t,\ACSBP s)\) for all spinors \((s,t)\in S=\SpinorBundleMEven\oplus\SpinorBundleMOdd\).
Hence, \(\SpinorBundleM\) may alternatively be viewed as a quaternionic line bundle.
This may not come as a big surprise for \(\CliffordAlgebraM\simeq \quaternions = \mathbb{R}\oplus\mathbb{R}^3\).
When viewed as complex vector bundles of rank two, the three complex spinor bundles \((\SpinorBundleM, \ACSBMI)\), \((\SpinorBundleM, \ACSBMJ)\) and \((\SpinorBundleM, \ACSBMK)\) are isomorphic and may be identified with \(\SpinorBundlePC = \ComplexHalfSpinorBundle\oplus\overline{\ComplexHalfSpinorBundle}\), whereby $\operatorname{Cl}(M,\pm g)\otimes\mathbb{C}\simeq_\mathbb{C}\operatorname{End}(\SpinorBundlePC)$.

Let us take a closer look at the identification of \((\SpinorBundleM, \ACSBMI)\) with \(\SpinorBundlePC\).
The spinor \((s, t)\in \SpinorBundleM=\SpinorBundleMEven\oplus\SpinorBundleMOdd\) is identified with \(s\otimes1 + t\otimes i\in \SpinorBundlePC=\SpinorBundleP\otimes\mathbb{C}\).
In particular \(\ACSBMI\) is identified with \(\Id_\SpinorBundleP\otimes i\).
Hence Equation~\eqref{eq:RepresentationM} can be rewritten as \(\gamma(X) = \gamma^+(X)\otimes i\), that is, the Clifford multiplication by~\(X\) on \(\SpinorBundleM\) differs from the Clifford multiplication by \(X\) on \(\SpinorBundleP\) by a factor of \(i\).
In this way any representation of \(\CliffordBundleP\) on \(\SpinorBundleP\) yields a purely imaginary representation of \(\CliffordBundleM\) on~\(\SpinorBundlePC\).
Furthermore, we obtain the following identifications of Dirac-operators:
\begin{equation}
	\DiracOperatorM = \DiracOperatorP\otimes i = i\DiracOperatorPC.
\end{equation}

We now derive a convenient local expression for the Dirac operator.
Let us first assume that \((M, g)\) is the Euclidean space with standard coordinates \(x\) and \(y\).
The holomorphic tangent bundle of \(M\) is then spanned by \(\partial_z = \frac12\left(\partial_x - i\partial_y\right)\).
The spinor bundle \((\SpinorBundleM, \ACSBMI)=\ComplexHalfSpinorBundle\oplus\overline{\ComplexHalfSpinorBundle}\) possesses a complex base \(s\), \(\overline{s}\) such that \(s\in\ComplexHalfSpinorBundle\), \(\overline{s}\) is the complex conjugate of \(s\) and \(s\otimes s =\dd{z}\).
With respect to this basis the Clifford multiplication of \(\CliffordBundleM\) on \((\SpinorBundleM, \ACSBMI)\) is represented by
\begin{align}
	\gamma(\partial_x) &=
	\begin{pmatrix}
		0 & 1 \\
		-1 & 0 \\
	\end{pmatrix},
	& \gamma\left(\partial_y\right) &=
	\begin{pmatrix}
		0 & -i \\
		-i & 0 \\
	\end{pmatrix}.
\end{align}
Hence the Euclidean Dirac-operator is given by
\begin{equation}
	\DiracOperatorM = 2
	\begin{pmatrix}
		0 & \partial_z \\
		-\partial_{\overline{z}} & 0 \\
	\end{pmatrix},
\end{equation}
that is, by the standard Cauchy--Riemann operators.
The general, non-Euclidean Dirac-operator differs from the Euclidean one by a rescaling and zero-order terms.
In particular, this means that the regularity theory developed for Cauchy--Riemann equations applies.

\subsection{}
In this paragraph we introduce the “super partner” of the metric, called gravitino.
\begin{Def}
	A \emph{gravitino} is a smooth section of the bundle \(\SpinorBundleM\otimes TM\).
\end{Def}
\begin{rmk}
	Sometimes in the literature, e.g.~\cite{jost2014super}, a gravitino is defined as a section of the bundle $\SpinorBundleM\otimes T^*M$, but here we use the Riemannian metric $g$ to identify $T^*M$ with $TM$, for later convenience.
\end{rmk}

The Clifford multiplication gives a surjective map
\begin{equation}
	\begin{split}
		\gamma\colon \SpinorBundleM\otimes TM&\to\SpinorBundleM \\
		s\otimes v&\mapsto v\cdot s
	\end{split}
\end{equation}
and has a canonical right-inverse that is given with respect to a local \(g\)-orthonormal base \(\{e_\alpha\}\) of \(TM\) by
\begin{equation}
	\begin{split}
		\sigma\colon\SpinorBundleM&\to \SpinorBundleM\otimes TM \\
		s&\mapsto -\frac12\delta^{\alpha\beta}e_\alpha\cdot s\otimes e_\beta.
	\end{split}
\end{equation}
Consequently the bundle \(\SpinorBundleM\otimes TM\) has an orthogonal direct sum decomposition \(\SpinorBundleM\otimes TM\cong\SpinorBundleM\oplus\ker\gamma\) and the maps \(P=\sigma\circ\gamma\) and \(Q=1-P\) are projection operators on \(\SpinorBundleM\) and \(\ker \gamma\) respectively.
With respect to the \(g\)-orthonormal frame \(\{e_\alpha\}\) the gravitino $\chi$ can locally be expressed as \(\chi = \chi^\alpha\otimes e_\alpha\) with $\chi^\alpha \in \Gamma_{loc}(\SpinorBundleM)$.
The projection operators \(P\) and \(Q\) are given by
\begin{align}
	P\chi &= -\frac{1}{2} e_\beta \cdot e_\alpha \cdot \chi^\alpha\otimes e_\beta, &
	Q\chi &= -\frac{1}{2} e_\alpha \cdot e_\beta \cdot \chi^\alpha \otimes e_\beta.
\end{align}

Later we will mostly be concerned with the sections of \(\ker\gamma\), because only \(Q\chi\) appears in the action functional.
Notice that \(\ker\gamma\) can be identified with \((\SpinorBundleM, \ACSBP\oplus\ACSBP)\otimes_{\mathbb{C}} TM\) because gravitinos of the form \(s\otimes\ACM v - (\ACSBP\oplus\ACSBP)s \otimes v\) span \(\ker\gamma\).
Using the almost complex structure \(\ACSBP\oplus\ACSBP\) on \(S\) and \(T^*M=W\otimes_\mathbb{C} W\) we obtain the following decomposition
\begin{equation}
	\begin{split}
		\SpinorBundleM\otimes TM &= \left(W\oplus W\right)\oplus \left(\left(W\oplus W\right) \otimes_\mathbb{C}\left(W^*\otimes_\mathbb{C} W^*\right)\right) \\
			&=W\oplus W\oplus \left(W\otimes_\mathbb{C} W^*\otimes_\mathbb{C} W^*\right)\oplus \left(W\otimes_\mathbb{C} W^*\otimes_{\mathbb{C}} W^*\right)
	\end{split}
\end{equation}
This is the decomposition of \(\SpinorBundleM\otimes TM\) into irreducible representations of \(\SpinGroup\).
Up to a metric identification, the bundle \(\SpinorBundleM\otimes TM\) decomposes into two representations of type \(\frac12\) and two of type \(\frac32\).
The operator \(Q\) projects onto the \(\frac32\)-parts.

\subsection{}
We recall the definition of the field \(\phi\) and its super partner \(\psi\), see~\cite{chen2006dirac}.
Let $(N,h)$ be a Riemannian manifold, with Levi-Civita connection $\na^N\equiv \na^{TN}$.
Consider a smooth map $\phi\colon M \to N$ with tangent map $T\phi\colon TM\to TN$.
It induces a pullback bundle~$\phi^*TN$ over~$M$.
Equip the tensor product bundle $S\otimes \phi^*TN$ with the induced metric and connection.
More precisely, let $\{y^i\}$ be local coordinates on N, so that $\{\phi^*(\frac{\p}{\p y^i})\}$ is a local frame of $\phi^* TN$.
Then the local sections, which will be referred to as ``(local) vector spinors'', can be written as $\psi=\psi^j\otimes \phi^*(\frac{\p}{\p y^j})$, $\varphi=\varphi^k\otimes \phi^*(\frac{\p}{\p y^k})$.
The induced metric and connection can be expressed by
\begin{gather}
	\langle \psi, \varphi \rangle_{S\otimes\phi^* TM}
	=\langle \psi^j,\varphi^k\rangle_{S} \cdot\big\langle \phi^*\frac{\p}{\p y^j}, \phi^*\frac{\p}{\p y^k}\big\rangle_{\phi^*TN}, \\
	\na^{S\otimes\phi^*TN}_X \psi
	= \na^S_X \psi^j \otimes\phi^*(\frac{\p}{\p y^j}) + \psi^j \otimes \na^{\phi^* TN}_X \phi^*(\frac{\p}{\p y^j}),
\end{gather}
where $\na^{\phi^* TN}_X \phi^*(\frac{\p}{\p y^j})=\phi^*(\na^{TN}_{T\phi(X)} \frac{\p}{\p y^j})$, for any $X\in TM$.
The twisted spin Dirac operator~$\D$ on $S\otimes\phi^*TN$ is defined as follows:
In a local \(g\)-orthonormal frame \(e_\alpha\) as above,
\begin{equation}
	\begin{split}
		\D\psi\coloneqq e_\al \cdot \na^{S\otimes \phi^*TN}_{e_\al} \psi
		&= e_\al \cdot \na^{S}_{e_\al} \psi^j \otimes \phi^*(\frac{\p}{\p y^j})
			+ e_\al \cdot \psi^j \otimes \na^{\phi^* TN}_{e_\al} \phi^*(\frac{\p}{\p y^j}) \\
		&= \pd \psi^j \otimes \phi^*(\frac{\p}{\p y^j})
			+ e_\al \cdot \psi^j \otimes \phi^*(\na^{ TN}_{T\phi e_\alpha} \frac{\p}{\p y^j}).
	\end{split}
\end{equation}
Similarly to the spin Dirac operator \(\pd\) the twisted spin Dirac operator \(\D\) is essentially self-adjoint with respect to the scalar product in \(L^2(S\otimes\phi^*TN)\).

\section{The Action Functional}%
\label{Sec:AF}
We want to consider the following action functional:
\begin{equation}
\label{eq:AF}
	\begin{split}
		\A(\phi, \psi;g, \chi)&\coloneqq \int_M |\dd \phi|_{T^*M\otimes \phi^*TN}^2
			+ \langle \psi, \D \psi \rangle_{S\otimes \phi^*TN} \\
			&\qquad -4\langle (\mathds{1}\otimes\phi_*)(Q\chi), \psi \rangle_{S\otimes\phi^*TN}
			-|Q\chi|^2_{S\otimes TM} |\psi|^2_{S\otimes \phi^*TN}
			-\frac{1}{6} \Rm^{\phi^*TN}(\psi) \dd vol_g,
	\end{split}
\end{equation}
where the last curvature term is locally defined by
\[-\frac{1}{6}\Rm^{\phi^*TN}(\psi)
=-\frac{1}{6}\Rm^{\phi^*TN}_{ijkl}\langle \psi^i, \psi^k \rangle_S \langle \psi^j,\psi^l \rangle_S. \]
Notice that we use the following conventions for the curvature tensor:
\begin{equation}
	\Rm^{TN}_{ijkl} = \left<\Rm^{TN}\left(\frac{\partial}{\partial y^k}, \frac{\partial}{\partial y^l}\right) \frac{\partial}{\partial y^j}, \frac{\partial}{\partial y^i}\right>
= \left<\nabla_{\frac{\partial}{\partial y^k}}\nabla_{\frac{\partial}{\partial y^l}}\frac{\partial}{\partial y^j} - \nabla_{\frac{\partial}{\partial y^l}}\nabla_{\frac{\partial}{\partial y^k}}\frac{\partial}{\partial y^j}, \frac{\partial}{\partial y^i}\right>
\end{equation}
We will abbreviate $\Rm^{\phi^*TN}$ as $\Rm^N$.
Hence, the curvature term can be written as
\begin{equation}
	\begin{split}
		\Rm^N(\psi)&=\Rm^N_{ijkl} \langle \psi^i, \psi^k \rangle_S \langle \psi^j, \psi^l \rangle_S
		=\Rm^N_{ijkl} \langle \psi^k, \psi^i \rangle_S \langle \psi^l, \psi^j \rangle_S \\
		&=\langle \Rm^N(\frac{\p}{\p y^k},\frac{\p}{\p y^l})\frac{\p}{\p y^j},\frac{\p}{\p y^i} \rangle_{TN}
			\langle \psi^k, \psi^i \rangle_S \langle \psi^l, \psi^j \rangle_S \\
		&=\big\langle \langle \psi^l, \psi^j \rangle_S \psi^k
			\otimes \phi^*(\Rm^N(\frac{\p}{\p y^k},\frac{\p}{\p y^l})\frac{\p}{\p y^j}),
			\psi^i\otimes \phi^*(\frac{\p}{\p y^i}) \big\rangle_{S\otimes \phi^*TN}.
	\end{split}
\end{equation}
So if we set
\begin{equation}%
\label{def-SR}
	SR(\psi)\coloneqq\langle \psi^l, \psi^j \rangle_S \psi^k
	\otimes \phi^*(\Rm^N(\frac{\p}{\p y^k},\frac{\p}{\p y^l})\frac{\p}{\p y^j}),
\end{equation}
then
\[ \Rm^N(\psi)=\langle SR(\psi), \psi \rangle_{S\otimes\phi^*TN}. \]

Note that since $P$ and $Q$ give an orthogonal decomposition,
\[ |Q\chi|^2_{S\otimes TM}= \langle \chi, Q\chi \rangle. \]
This formula is convenient when expressing the terms locally.

\begin{rmk}
	In order to obtain a real-valued action functional we work here with the real spinor bundle \(\SpinorBundleM\) and the real scalar product \(\SpinorMetricM=\left<\cdot, \cdot\right>_\SpinorBundleM\).
	Alternatively we might also work with the complex spinor bundle \(\SpinorBundlePC\) and the hermitian form \(\SpinorHermitianFormPC\).
	We recall that the hermitian form \(\SpinorHermitianFormM\) on \((\SpinorBundleM, \ACSBMI)\) induced by \(\SpinorMetricM\) can be written as
	\begin{equation}
		2\SpinorHermitianFormM\left(s, t\right) = \SpinorMetricM\left(s, t\right) - i\SpinorMetricM\left(\ACSBMI s, t\right)
	\end{equation}
	and coincides with \(\SpinorHermitianFormPC\) under the complex linear isomorphism \(\SpinorBundlePC\simeq(\SpinorBundleM, \ACSBMI)\).
	All summands in~\eqref{eq:AF} except the third one are symmetric in the spinors and would consequently be real.
	For those terms the approach here and in~\cite{chen2006dirac} coincide.
	For the third term one could use equally the real part of
	\begin{equation}
		-8\SpinorHermitianFormPC\otimes\phi^*h\left((\mathds{1}\otimes\phi_*)(Q\chi), \psi\right).
	\end{equation}
	We will refrain from using that expression later on.
\end{rmk}

The functional $\A(\phi,\psi;g,\chi)$ has rich symmetries.
It is invariant under generalized conformal transformations of the metric in the sense that
\begin{equation}
	\A(\phi,e^{-u}\psi;e^{2u}g,e^{-2u}\chi) =\A(\phi,\psi;g,\chi)
\end{equation}
where $u\in C^\infty(M)$.
To verify the conformal invariance we use the rescaling of the spinor metric~\(\SpinorMetricM\) by \(e^u\SpinorMetricM\) and that \(\D^{e^{2u}g}e^{-u}\psi = e^{-2u}\D^g\psi\), see also~\cite[Proposition~1.3.10]{ginoux2009}.
Here \(\D^g\) denotes the Dirac operator defined with respect to the metric \(g\).
Moreover, the functional stays invariant under super Weyl transformations:
\begin{equation}
	\A(\phi,\psi;g,\chi+\chi')=\A(\phi,\psi;g,\chi)
\end{equation}
with $Q\chi'=0$.
This follows directly from the fact that the action functional only involves \(Q\chi\) and not \(P\chi\).
$\A$ is also \(\SpinGroup\)-gauge-invariant, in particular under the following $\mathbb{Z}_2$-action on the spinor bundle $S$:
\begin{equation}
	\A(\phi,\psi;g,\chi)=\A(\phi,-\psi;g,-\chi).
\end{equation}
These symmetries will be naturally inherited by its critical points.
They are useful when dealing with the solution space of the Euler--Lagrange equations.
A detailed discussion of the symmetries of \(\A\) and the corresponding conservation laws can be found in~\cite{jost2017symmetries}.

As already mentioned in the introduction the functional~\eqref{eq:AF} is essentially the action functional of the two-dimensional nonlinear supersymmetric sigma model, see~\cite{brink1976locally, deser1976complete, jost2014super}.
In contrast to what is discussed there, we deal with commuting spinors.
For that matter the action functional~\eqref{eq:AF} does in general not possess supersymmetry, except in special cases, see~\cite{jost2017symmetries}.
Furthermore, a term which vanishes identically at critical points is omitted here.

\section{Euler--Lagrange Equations}

\subsection{}
Now we derive the Euler--Lagrange equations for $\A$.
Fix $(g,\chi)$ and vary $(\phi, \psi)$ via $(\Phi, \Psi)$ with variational fields $(\xi, \eta)$.
Here
\begin{align}
	\xi &= \left.\frac{\partial}{\partial t}\Phi\right|_{t=o}, &
	\eta &= \left.\nabla_{\partial_t}^{S\otimes\Phi^*TN}\Psi\right|_{t=0}.
\end{align}
At a critical point, we have
\[ 0=\frac{\dd}{\dd t}\Big|_{t=0} \A(\Phi(t),\Psi(t);g, \chi) = \frac{\dd}{\dd t}\Big|_{t=0} (\textrm{I+II+III+IV+V}).\]
Here we denote by the roman numerals \(\textrm{I}, \dotsc, \textrm{V}\) the summands under the integral in the action functional \(\A\).
We calculate them term by term.
\begin{enumerate}
	\item As for harmonic maps,
		\[\frac{\dd}{\dd t} \textrm{I}
		= \frac{\dd}{\dd t} \int_M |d_x\Phi|^2= \int_M \langle -2\tau(\Phi), \Phi_* (\p_t) \rangle_{\Phi^*TN}, \]
		where $\tau(\Phi)$ is the tension field of $\Phi$ w.r.t.\ $M$. Hence,
		\[ \frac{\dd}{\dd t}\Big|_{t=0} \textrm{I}=\int_M \langle -2\tau(\phi), \xi \rangle_{\phi^*TN}.\]
	\item With
		\[ \na^{S\otimes\Phi^*TN}_{\p_t} \D \Psi
		= \D \na^{S\otimes\Phi^*TN}_{\p_t} \Psi+\Rm^{\Phi^*TN} (\Phi(\p_t), \Phi_* e_\al)e_\al \cdot \Psi,\]
		we get
		\begin{equation}
			\begin{split}
				\frac{\dd}{\dd t} \textrm{II}
				&=\frac{\dd}{\dd t} \int_M \langle \Psi, \D \Psi \rangle_{S\otimes \Phi^*TN}
				=\int_M \langle \na^{S\otimes\Phi^*TN}_{\p_t} \Psi, \D \Psi \rangle
					+ \langle \Psi, \na^{S\otimes\Phi^*TN}_{\p_t}\D \Psi \rangle \\
				&= \int_M \langle \na^{S\otimes\Phi^*TN}_{\p_t} \Psi, \D \Psi \rangle
					+ \langle \Psi, \D \na^{S\otimes\Phi^*TN}_{\p_t} \Psi
					+\Rm^{\Phi^*TN} (\Phi_*(\p_t), \Phi_* e_\al)e_\al \cdot \Psi \rangle \\
				&= \int_M \langle \na^{S\otimes\Phi^*TN}_{\p_t} \Psi, \D \Psi \rangle
					+ \langle \D \Psi,\na^{S\otimes\Phi^*TN}_{\p_t} \Psi\rangle
					+\langle \Psi,\Rm^{\Phi^*TN} (\Phi_*(\p_t), \Phi_* e_\al)e_\al \cdot \Psi \rangle \\
				&= \int_M 2 \langle \na^{S\otimes\Phi^*TN}_{\p_t} \Psi, \D \Psi \rangle
					+\langle \Rm^{\Phi^*TN} (\Psi, e_\al \cdot\Psi)\Phi_* e_\al ,\Phi_*(\p_t)\rangle.
			\end{split}
		\end{equation}
		Thus
		\[\frac{\dd}{\dd t}\Big|_{t=0} \textrm{II}
			=\int_M 2\langle \eta, \D \psi \rangle
			+ \langle \Rm^{\phi^*TN}(\psi, e_\al \cdot \psi)\phi_* e_\al, \xi \rangle. \]
	\item Under a local orthonormal frame $\{ e_\al\}$,
		\[ -4\langle (\mathds{1}\otimes\Phi_*)(Q\chi), \Psi \rangle_{S\otimes\Phi^*TN}
			=2\langle e_\al \cdot e_\be \cdot \chi^\al \otimes \Phi_* e_\be, \Psi \rangle_{S\otimes \Phi^*TN}. \]
		Then
		\begin{equation}
			\begin{split}
				& \frac{\dd}{\dd t} \textrm{III}
				=\frac{\dd}{\dd t} \int_M 2\langle e_\al \cdot e_\be \cdot \chi^\al \otimes \Phi_* e_\be,
				\Psi \rangle_{S\otimes \Phi^*TN} \\
				&\;=\int_M 2\langle \na^{S\otimes\Phi^* TN}_{\p_t}(e_\al\cdot e_\be\cdot\chi^\al\otimes\Phi_* e_\be), \Psi\rangle
				+2\langle e_\al\cdot e_\be\cdot\chi^\al\otimes\Phi_* e_\be,\na^{S\otimes\Phi^*TN}_{\p_t}\Psi \rangle,
			\end{split}
		\end{equation}
		where the first integrand can be rewritten as
		\begin{equation}
			\begin{split}
				2\langle \na^{S\otimes \Phi^* TN}_{\p_t}&(e_\al \cdot e_\be \cdot \chi^\al \otimes \Phi_* e_\be), \Psi \rangle
				= 2\langle e_\al \cdot e_\be \cdot \chi^\al \otimes \na^{\Phi^* TN}_{\p_t}\Phi_* e_\be, \Psi \rangle \\
				&= 2\langle e_\al \cdot e_\be \cdot \chi^\al \otimes \na^{\Phi^* TN}_{e_\be}\Phi_* \p_t, \Psi \rangle \\
				&= 2 e_\be \langle e_\al \cdot e_\be \cdot \chi^\al \otimes \Phi_* \p_t, \Psi \rangle -2\langle \na^S_{e_\be}( e_\al \cdot e_\be \cdot \chi^\al) \otimes\Phi_* \p_t, \Psi \rangle \\
				&\quad -2\langle e_\al \cdot e_\be \cdot \chi^\al \otimes\Phi_* \p_t, \na^{S\otimes \Phi^*TN}_{e_\be} \Psi \rangle.
			\end{split}
		\end{equation}
		The first summand vanishes after integration on the closed manifold~$M$ since it is a divergence of some vector field.
		Therefore
		\begin{equation}
			\begin{split}
				\frac{\dd}{\dd t}\Big|_{t=0} \textrm{III}
				=\int_M &-2\langle  \na^S_{e_\be}( e_\al \cdot e_\be \cdot \chi^\al) \otimes \xi, \psi \rangle
				-2\langle e_\al \cdot e_\be\cdot \chi^\al \otimes\xi,\na^{S\otimes \phi^*TN}_{e_\be}\psi \rangle \\
				& +2\langle e_\al \cdot e_\be \cdot \chi^\al \otimes \phi_* e_\be, \eta \rangle \\
				=\int_M &-2\big\langle(\langle \na^S_{e_\be}(e_\al \cdot e_\be \cdot \chi^\al), \psi \rangle_S
					+\langle e_\al\cdot e_\be\cdot\chi^\al,\na^{S\otimes\phi^*TN}_{e_\be}\psi\rangle_S),
					\xi \big\rangle_{\phi^*TN} \\
				& +2\langle e_\al \cdot e_\be \cdot \chi^\al \otimes \phi_* e_\be, \eta \rangle_{S\otimes \phi^*TN}.
			\end{split}
		\end{equation}
		Here, by abuse of notation we denote by \(\langle \na^S_{e_\be}(e_\al \cdot e_\be \cdot \chi^\al), \psi \rangle_S\), the section of \(\phi^*TN\) that arises by metric contraction of \(\psi\) by \(\na^S_{e_\be}(e_\al \cdot e_\be \cdot \chi^\al)\).
	\item Likewise we have
		\begin{equation}
			\begin{split}
				\frac{\dd}{\dd t} \textrm{IV}
				=& -\frac{\dd}{\dd t} \int_M |Q\chi|^2 \langle \Psi, \Psi \rangle_{S\otimes\Phi^*TN}  \\
				=& -\int_M |Q\chi|^2 (\langle \na^{S\otimes \Phi^* TN}_{\p_t} \Psi, \Psi \rangle
				+ \langle \Psi,\na^{S\otimes \Phi^* TN}_{\p_t} \Psi  \rangle) \\
				=& -\int_M 2|Q\chi|^2 \langle \Psi, \na^{S\otimes \Phi^* TN}_{\p_t} \Psi \rangle.
			\end{split}
		\end{equation}
		Thus,
		\[ \frac{\dd}{\dd t}\Big|_{t=0} \textrm{IV} = \int_M -2|Q\chi|^2 \langle \psi, \eta \rangle. \]
	\item In local coordinates, we compute
		\begin{equation}
			\begin{split}
				\frac{\dd}{\dd t} \textrm{V}
				&=\frac{\dd}{\dd t}\int_M-\frac{1}{6}\Phi^*\Rm^N_{ijkl}\langle\Psi^i,\Psi^k\rangle_S \langle\Psi^j,\Psi^l\rangle_S\\
				&= -\frac{1}{6}\int_M \p_t (\Phi^*\Rm^N_{ijkl}\langle \Psi^i, \Psi^k \rangle_S \langle \Psi^j, \Psi^l \rangle_S).
			\end{split}
		\end{equation}
		The integrand reads
		\begin{equation}
			\begin{split}
				\p_t (\Phi^*\Rm^N_{ijkl}&\langle \Psi^i, \Psi^k \rangle \langle \Psi^j, \Psi^l \rangle)  \\
				= &(\na^{\Phi^*TN}_{\p_t} \Phi^*\Rm^N_{ijkl})\langle\Psi^i,\Psi^k\rangle \langle\Psi^j,\Psi^l\rangle \\
				&   + \Phi^*\Rm^N_{ijkl}\langle \na^S_{\p_t}\Psi^i, \Psi^k \rangle \langle \Psi^j, \Psi^l \rangle
					+ \Phi^*\Rm^N_{ijkl}\langle \Psi^i, \na^S_{\p_t} \Psi^k \rangle \langle \Psi^j, \Psi^l \rangle  \\
				&   + \Phi^*\Rm^N_{ijkl}\langle \Psi^i, \Psi^k \rangle \langle \na^S_{\p_t} \Psi^j, \Psi^l \rangle
					+ \Phi^*\Rm^N_{ijkl}\langle \Psi^i, \Psi^k \rangle \langle \Psi^j, \na^S_{\p_t} \Psi^l \rangle  \\
				= & (\na^{\Phi^*TN}_{\p_t} \Phi^*\Rm^N_{ijkl})\langle\Psi^i,\Psi^k\rangle \langle\Psi^j,\Psi^l\rangle
					+4\Phi^*\Rm^N_{ijkl}\langle \na^S_{\p_t}\Psi^i,\Psi^k\rangle \langle\Psi^j,\Psi^l\rangle\\
				= &\Phi^*(\na^{TN}_{T\Phi(\p_t)}\Rm^N)_{ijkl}\langle\Psi^i,\Psi^k\rangle \langle\Psi^j, \Psi^l \rangle
					+4\langle \na^{S\otimes\Phi^*TN}_{\p_t} \Psi, SR(\Psi) \rangle  \\
				= &\big\langle \Phi^*(\na^{TN}\Rm^N)_{ijkl}\langle\Psi^i,\Psi^k\rangle \langle\Psi^j,\Psi^l\rangle, \Phi_*\p_t \big\rangle
					+4 \langle \na^{S\otimes\Phi^*TN}_{\p_t} \Psi, SR(\Psi) \rangle.  \\
			\end{split}
		\end{equation}
		We define $S\na R$ analogously to $SR$, that is,
		\begin{equation}\label{def-SNR}
			S\na R(\Psi)\coloneqq \Phi^*(\na^{TN} \Rm^N)_{ijkl}\langle \Psi^i, \Psi^k \rangle \langle \Psi^j, \Psi^l \rangle.
		\end{equation}
		Using the metric to identify it with the corresponding vector field, we get
		\begin{equation}
			\p_t (\Phi^*\Rm^N_{ijkl}\langle \Psi^i, \Psi^k \rangle \langle \Psi^j, \Psi^l \rangle)
			= \langle S\na R (\Psi), \Phi_* \p_t \rangle
			+ 4 \langle \na^{S\otimes\Phi^*TN}_{\p_t} \Psi, SR(\Psi) \rangle.
		\end{equation}
		Then,
		\[ \frac{\dd}{\dd t}\Big|_{t=0} \textrm{V}
		= -\frac{1}{6}\int_M \langle S\na R (\psi), \xi \rangle + 4\langle \eta, SR(\psi) \rangle.  \]
\end{enumerate}

From the preceding computations, we obtain
\begin{equation}
	\begin{split}
		0= \int_M &\big\langle -2\tau(\phi)+\Rm^N(\psi,e_\al\cdot\psi)\phi_* e_\al-\frac{1}{6}S\na R(\psi), \xi \big\rangle_{\phi^*TN} \\
		&+\big\langle -2(\langle \na^S_{e_\be}(e_\al \cdot e_\be \cdot \chi^\al), \psi \rangle_S
			+ \langle e_\al \cdot e_\be \cdot \chi^\al, \na^{S\otimes\phi^*TN}_{e_\be} \psi \rangle_S), \xi \big\rangle_{\phi^*TN}\\
		&+2 \langle \D \psi -|Q\chi|^2 \psi -\frac{1}{3}SR(\psi)
			+e_\al \cdot e_\be \cdot \chi^\al \otimes \phi_* e_\be,
			\eta \rangle_{S\otimes \phi^*TN}.
	\end{split}
\end{equation}
We can thus verify Theorem~\ref{thm:EL} which we restate here:
\ELThm*

\begin{Def}
	A pair $(\phi,\psi)\in \mathcal{X}^{1,2}_{1,4/3}(M,N)$ satisfying~\eqref{EL-eq} in the sense of distributions is a weak solution of the system.
\end{Def}

\subsection{}
We rewrite the Euler--Lagrange equations~\eqref{EL-eq} in terms of local coordinates on $N$.

Let $\{ y^i \} $ be a local coordinate system on $N$. Then $\{ \phi^*(\frac{\p}{\p y^i}) \}$ is a local frame for
the vector bundle $\phi^*TN$.
Then~\eqref{EL-eq} can be written as
\begin{equation}
	\begin{split}
		\tau(\phi)^i \phi^*(\frac{\p}{\p y^i})
		=& \frac{1}{2}\langle\psi^k, e_\al\cdot\psi^l\rangle
			\Rm^N\big(\frac{\p}{\p y^k},\frac{\p}{\p y^l}\big)\big(e_\al(\phi^j)\phi^*(\frac{\p}{\p y^j})\big)
			-\frac{1}{12}(\na \Rm^N)_{mjkl} \langle \psi^m,\psi^k \rangle \langle \psi^j, \psi^l\rangle  \\
		& \quad -\big( \langle \na^S_{e_\be} (e_\al \cdot e_\be \cdot \chi^\al), \psi^i \rangle
			+ \langle e_\al\cdot e_\be \cdot\chi^\al,\na^S_{e_\be}\psi^i\rangle \big)\phi^*(\frac{\p}{\p y^i}) \\
		& \quad-\langle e_\al \cdot e_\be\cdot \chi^\al,\psi^k \rangle \na^{\phi^*TN}_{e_\be} \phi^*(\frac{\p}{\p y^k})  \\
		=& \Big(
			\frac{1}{2}\langle \psi^k, e_\al \cdot \psi^l \rangle e_\al(\phi^j)\Rm^{i,N}_{\; jkl}
			-\frac{1}{12}(\na^i \Rm^N)_{mjkl} \langle \psi^m,\psi^k \rangle \langle \psi^j, \psi^l\rangle  \\
		&    \quad - e_\be( \langle e_\al \cdot e_\be \cdot \chi^\al,\psi^i \rangle)
			-\langle e_\al \cdot e_\be \cdot \chi^\al, \psi^k \rangle e_\be(\phi^j) \Gamma^{i,N}_{jk}
		\Big) \phi^*(\frac{\p}{\p y^i})
	\end{split}
\end{equation}
and
\begin{equation}
	\begin{split}
		\pd\psi^i \otimes &\phi^*(\frac{\p}{\p y^i})
		+ e_\al \cdot \psi^k \otimes e_\al(\phi^j)\Gamma^{i,N}_{jk} \phi^*(\frac{\p}{\p y^i})
		\\
		&= |Q\chi|^2 \psi^i \otimes \phi^*(\frac{\p}{\p y^i})
		+\frac{1}{3} \langle \psi^l, \psi^j \rangle \psi^k\otimes \Rm^{i,N}_{jkl} \phi^*(\frac{\p}{\p y^i})
		-e_\al \cdot e_\be \cdot \chi^\al \otimes e_\be(\phi^i) \phi^*(\frac{\p}{\p y^i}).
	\end{split}
\end{equation}
Since the curvature of $M$ does not appear in those formulas, we may omit the upper index
$N$ for the curvature terms, and we will label it again whenever needed.

We may introduce local coordinates on $M$ such that a conformal transformation brings the metric into the following form
\[ g= (\dd x^1)^2 + (\dd x^2)^2,\]
and then $\{e_\al\equiv \frac{\p}{\p x^\al}\}$ is a local orthonormal frame.
We define the vector fields $V^j$ on $M$, $j=1,\dotsc, n$, via
\[ \langle V^j,W \rangle_{TM}= \langle e_\al \cdot W \cdot \chi^\al, \psi^j \rangle_{S} \]
for any vector field $W$ on $M$.
Thus,
\[ V^j= V^{j,\be}e_\be = \langle e_\al \cdot e_\be \cdot \chi^\al, \psi^j \rangle e_\be.\]
In particular, noting that $\na_{e_\al} e_\be=0$, we have
\[ \diverg  V^j= e_\be(V^{j,\be})= e_\be \langle e_\al \cdot e_\be \cdot \chi^\al, \psi^j \rangle.  \]
and
\[\langle e_\al \cdot e_\be \cdot \chi^\al, \psi^k \rangle e_\be(\phi^j) \Gamma^{i,N}_{jk}
=V^{k,\be} e_{\be}(\phi^j) \Gamma^i_{jk}= \Gamma^i_{jk}V^k(\phi^j)=\Gamma^i_{jk}\langle V^k,\na\phi^j\rangle_{TM}.\]
Thus, in those local coordinates the Euler--Lagrange equations become
\begin{equation}\label{local form of EL-eq on N}
	\begin{split}
		\tau^i(\phi)=& \frac{1}{2}\langle \psi^k, e_\al \cdot \psi^l \rangle e_\al(\phi^j)\Rm^{i}_{\; jkl}
		-\frac{1}{12}(\na^i \Rm)_{mjkl} \langle \psi^m,\psi^k \rangle \langle \psi^j, \psi^l\rangle  \\
		& -\diverg V^i -\Gamma^{i}_{jk}\langle V^k,\na\phi^j \rangle, \\
		\pd \psi^i  =& -\Gamma^i_{jk}\na \phi^j\cdot\psi^k+|Q\chi|^2\psi^i
		+\frac{1}{3}R^i_{\;jkl}\langle\psi^l,\psi^j\rangle\psi^k
		-e_\al\cdot \na \phi^i \cdot \chi^\al,
	\end{split}
\end{equation}
for $1\le i \le n$. One sees that the right hand side of the first equation lies in $L^1$ while
that of the second equation lies in $L^{4/3}$.
This shows that the Euler--Lagrange equations are critical for the Sobolev elliptic theory.
Thus, the regularity of weak solutions is a subtle issue.

\subsection{}
To get the regularity of weak solutions, we embed $N$ isometrically into some Euclidean space.
In order to see what happens to the various fields involved, we start with a general consideration.
Let $(N',h')$ be another Riemannian manifold and $f\colon N\to N'$ a smooth immersion.
We get a composition
\[\phi'\equiv f\circ\phi\colon M\to N\to N',\]
and induced maps of vector bundles which fit into the following commutative diagram
\begin{diag}
	\matrix[mat, column sep=small](m)
			{
					& TM &   &          & (f\circ\phi)^*TN'&     & f^*TN'&     &   TN' \\
				TM &    &   & \phi^*TN &                  & TN  &       &     &      \\
					&    &   &          &                  &     &       &     &      \\
					&    & M &          &                  &     &  N    &     &   N'  \\
				} ;
	\path[pf]
			(m-1-2) edge node{\( (f\circ\phi)_* \)} (m-1-5)
				edge [dashed](m-4-3)
			(m-1-2) edge [commutative diagrams/equal](m-2-1)
			(m-1-5) edge node[auto]{\(\Hat{\Hat{\phi}}\)} (m-1-7)
							edge [dashed](m-4-3)
			(m-1-7) edge node[auto]{\(\Hat{f} \)} (m-1-9)
							edge [dashed](m-4-7)
			(m-1-9) edge (m-4-9)
			(m-2-1) edge [near end]node[auto]{\(\quad\phi_* \)} (m-2-4)
							edge (m-4-3)
			(m-2-4) edge [commutative diagrams/crossing over]node[auto]{\(\Hat{\phi} \)} (m-2-6)
							edge (m-4-3)
							edge node[auto]{\(\Hat{\phi}^*(f_*)\)} (m-1-5)
			(m-2-6) edge node[auto]{\( f_* \)} (m-1-7)
							edge (m-4-7)
			(m-4-3) edge node[auto]{\( \phi \)}  (m-4-7)
			(m-4-7) edge node[auto]{\( f \)}  (m-4-9)
						;
\end{diag}
Note that $T\phi=\Hat{\phi}\circ \phi_*$, etc.
Let $A$ be the second fundamental form of $f$, i.e., $A(X,Y)=(\na_X \dd f)(Y)$ for any $X,Y\in\Gamma(TN)$.
Then the tension fields of $\phi$ and $\phi'$ are related  by
\begin{equation}
	\tau(\phi')= \Hat{\phi}^*(f_*)(\tau(\phi))+A(\phi)\big(T\phi( e_\al),T\phi(e_\al)\big).
\end{equation}

Now let $(N',h')=(\R^K,\delta)$ be a Euclidean space with standard global coordinate functions $(u^a)_{a=1,\dotsc,K}$, and let $f\colon(N,h)\to (\R^K,\delta)$ be an isometric embedding.
Then the second fundamental form $A$ is perpendicular to $N$ in the sense that, for any $X,Y\in\Gamma(TN)$, extended locally to $\R^K$ and still denoted by $X,Y$ respectively, the following orthogonal decomposition holds:
\begin{equation}
	\na^e_X Y=\na^N_X Y+ A(X,Y)\in TN\oplus T^\bot N=f^*T\R^K,
\end{equation}
where $\na^e$ denotes the flat connection on Euclidean space.
See~\cite{bai2004anintroductiontoRG, jost2008riemannian}.
Moreover, for any normal vector field
$\xi\in \Gamma(T^\bot N)$,
\begin{equation}
	\langle\xi, A(X,Y)\rangle=\langle\xi,\na^e_X Y\rangle=-\langle\na^e_X\xi,Y\rangle=\langle P(\xi;X),Y\rangle
\end{equation}
where $P(\xi;X)=-(\na^e_X \xi)^\top$ is the shape operator of $N$.

As in~\cite{zhu2009regularity} and~\cite{chen2011boundary}, we take a local orthonormal frame
$\{\nu_l|l=n+1,\dotsc,K\}$ of $T^\bot N$. (These can be smoothly extended to a tubular neighborhood of $N$,
and thus be defined in an open subset of $\R^K$).
Then
\begin{equation}
	A(X,Y)=\sum_l \langle A(X,Y),\nu_l\rangle\nu_l=-\sum_{l} \langle Y,\na^e_X \nu_l\rangle \nu_l.
\end{equation}
In terms of the global frame $\{\frac{\p}{\p u^a}\}$ we write the vector fields $X,Y,Z$ tangent to the submanifold $N$
as
\begin{equation}
	X=X^a \frac{\p}{\p u^a}, \quad Y=Y^b \frac{\p}{\p u^b}, \quad Z=Z^c \frac{\p}{\p u^c}.
\end{equation}
Then
\begin{equation}
	\begin{split}
		A(X,Y)=\;&\sum_l -\langle Y^b\frac{\p}{\p u^b},\na^e_{X^a\frac{\p}{\p u^a}}\nu_l\rangle\nu_l
		=-\sum_{l,b}X^a Y^b \frac{\p \nu_l^b}{\p u^a}\nu_l,\\
		P(A(X,Y);Z)=\;& -(\na^e_Z A(X,Y))^\top=\sum_{l,b} Z^c X^a Y^b \frac{\p\nu_l^b}{\p u^a}(\frac{\p \nu_l}{\p u^c})^\top.
	\end{split}
\end{equation}
Since $A$ is symmetric: $A(X,Y)=A(Y,X)$, we have
\begin{equation}\label{symmetry of A}
	A(X,Y)=-\sum_{l,b}X^a Y^b \frac{\p \nu_l^b}{\p u^a}\nu_l=-\sum_{l,b}X^b Y^a \frac{\p \nu_l^b}{\p u^a}\nu_l,
\end{equation}
\begin{equation}\label{symmetry of P}
	P(A(X,Y);Z)=\sum_{l,b} Z^c X^a Y^b \frac{\p\nu_l^b}{\p u^a}(\frac{\p \nu_l}{\p u^c})^\top
	=\sum_{l,b} Z^c X^b Y^a \frac{\p\nu_l^b}{\p u^a}(\frac{\p \nu_l}{\p u^c})^\top.
\end{equation}
We recall here the Gauss equation for $X,Y,Z,W\in \Gamma(TN)$:
\begin{equation}
	\begin{split}
		\left<\Rm\left(X, Y\right)Z, W\right>
		&= \left<A\left(X, W\right), A\left(Y,Z\right)\right> - \left<A\left(X, Z\right), A\left(Y, W\right)\right> \\
		&= \left<P\left(A\left(Y, Z\right); X\right), W\right> - \left<P\left(A\left(X, Z\right); Y\right), W\right>
	\end{split}
\end{equation}
Since this holds for all $W\in \Gamma(TN)$, we have
\begin{equation}%
\label{Rm-P}
	\Rm\left(X, Y\right)Z = P\left(A\left(Y, Z\right); X\right) - P\left(A\left(X, Z\right); Y\right).
\end{equation}

We will denote the induced map on the tensor product bundles by
\begin{equation}
	f_\#\equiv\mathds{1}\otimes \Hat{\phi}^*(f_*)\colon S\otimes\phi^*TN\to S\otimes\phi'^*TN'.
\end{equation}
Then $\psi'\equiv f_\#(\psi)$ is a section of the latter bundle, i.e., a spinor field along the map $\phi'$.
In local coordinates,
\begin{equation}
	\psi=\psi^i\otimes\phi^*(\frac{\p}{\p y^i}), \quad \psi'=\psi'^a\otimes \phi'^*(\frac{\p}{\p u^a}),
\end{equation}
where
\begin{equation}\label{psi'-psi}
	\psi'^a(x)=\frac{\p u^a}{\p y^i}(\phi(x))\psi^i(x).
\end{equation}
Moreover, the Dirac terms corresponding to $\phi$ and $\phi'$ are related via (see~\cite{chen2011boundary})
\begin{equation}\label{D'-D}
	\D'\psi'=f_\# \D\psi+\mathcal{A}(\phi_* e_\al, e_\al\cdot\psi),
\end{equation}
where
\begin{equation}
	\mathcal{A}(\phi_* e_\al, e_\al\cdot\psi)
	\equiv e_\al\cdot\psi^i\otimes \phi^*\big(A(T\phi (e_\al),\frac{\p}{\p y^i})\big).
\end{equation}

\subsection{}
We are now ready to write the Euler--Lagrange equations in terms of $(\phi',\psi')$.


Apply $f_\#$ to $\D \psi$ and use~\eqref{D'-D}:
\begin{equation}
	\begin{split}
		\D'\psi'-\mathcal{A}(\phi_* e_\al,e_\al\cdot\psi)=|Q\chi|^2\psi'+\frac{1}{3}f_\# (SR(\psi))
		+2(\mathds{1}\otimes\phi'_*)Q\chi.
	\end{split}
\end{equation}
We compute the following terms:
\begin{itemize}
	\item Note that
		\begin{equation}
			\begin{split}
				Tf \left(\frac{\p}{\p y^i}\right) = \frac{\p f^a}{\p y^i}\frac{\p}{\p u^a}=\frac{\p u^a}{\p y^i}(\phi)\frac{\p}{\p u^a}
			\end{split}
		\end{equation}
		and
		\begin{equation}
			T\phi'(e_\al)=\frac{\p \phi^i}{\p x^\al}Tf\left(\frac{\p}{\p y^i}\right)
			=\frac{\p\phi^i}{\p x^\al}\frac{\p f^a}{\p y^i}\frac{\p}{\p u^a}
			=\frac{\p \phi'^a}{\p x^\al}\frac{\p}{\p u^a}.
		\end{equation}
		Using~\eqref{psi'-psi} and the expression for $A$, we have
		\begin{equation}
			\begin{split}
				\mathcal{A}(\phi_* e_\al, e_\al\cdot\psi)
				=\;& e_\al\cdot\psi^i\otimes \phi^*(A(T\phi(e_\al),\frac{\p}{\p y^i})) \\
				=\;& -e_\al\cdot\psi^i\otimes
				\sum_{l,b}\frac{\p\phi'^a}{\p x^\al} \frac{\p u^b}{\p y^i}(\phi)
				\frac{\p \nu_l^b}{\p u^a}(\phi') \phi'^* \nu_l \\
				=\;& -\sum_{l,b} \frac{\p \phi'^a}{\p x^\al}e_\al \cdot \frac{\p u^b}{\p y^i}\psi^i
				\otimes \frac{\p \nu_l^b}{\p u^a} \nu_l^c(\phi')\phi'^*(\frac{\p}{\p u^c}) \\
				=\;& -\sum_{l,b} \na\phi'^a \cdot \psi'^b
				\otimes \frac{\p \nu_l^b}{\p u^a} \nu_l^c(\phi')\phi'^*(\frac{\p}{\p u^c}).
			\end{split}
		\end{equation}
	\item Recalling~\eqref{def-SR} and~\eqref{Rm-P},
		\begin{equation}
			\begin{split}
				f_\# SR(\psi)
				=\;&f_\#\big(\langle\psi^l,\psi^j\rangle\psi^k
					\otimes\Rm(\frac{\p}{\p y^k},\frac{\p}{\p y^l})\frac{\p}{\p y^j}\big) \\
				=\;& f_\#\big\{\langle\psi^l,\psi^j\rangle\psi^k
					\otimes \big(P(A(\frac{\p}{\p y^j},\frac{\p}{\p y^l});\frac{\p}{\p y^k})
					-P(A(\frac{\p}{\p y^k},\frac{\p}{\p y^j});\frac{\p}{\p y^l})\big)\big\}\\
				=\;& \langle\psi^l,\psi^j\rangle\psi^k
					\otimes\frac{\p u^a}{\p y^j}\frac{\p u^b}{\p y^l}\frac{\p u^c}{\p y^k}
					\frac{\p \nu_l^b}{\p u^a}(\frac{\p \nu_l}{\p u^c})^\top(\phi')   \\
				&\quad -\langle\psi^l,\psi^j\rangle\psi^k
					\otimes\frac{\p u^a}{\p y^k}\frac{\p u^b}{\p y^j}\frac{\p u^c}{\p y^l}
					\frac{\p \nu_l^b}{\p u^a}(\frac{\p \nu_l}{\p u^c})^\top(\phi')  \\
				=\;& \big(\langle\psi'^b,\psi'^a\rangle\psi'^c
					-\langle\psi'^c,\psi'^b \rangle\psi'^a\big)\otimes
					\frac{\p \nu_l^b}{\p u^a}(\frac{\p \nu_l}{\p u^c})^{\top,d}\phi'^*(\frac{\p}{\p u^d}).
			\end{split}
		\end{equation}
	\item For the last term:
		\begin{equation}
			\begin{split}
				2(\mathds{1}\otimes\phi'_*)Q\chi
				=\;& -e_\al\cdot e_\be\cdot \chi^\al \otimes\phi'_* e_\be
				=-e_\al\cdot e_\be\cdot\chi^\al\otimes \frac{\p \phi'^a}{\p x^\be}\phi'^*(\frac{\p}{\p u^a})\\
				=\;& -e_\al\cdot \na \phi'^a\cdot\chi^\al \otimes \phi'^*(\frac{\p}{\p u^a}).
			\end{split}
		\end{equation}
\end{itemize}
We thus obtain the equation for $\psi'$:
\begin{equation}%
\label{eq-psi'}
	\begin{split}
		\pd\psi'^a\otimes\phi'^*(\frac{\p}{\p u^a})
		=\;& -\sum_{l,b}\na\phi'^d\cdot\psi'^b\otimes\frac{\p\nu_l^b}{\p u^d}\nu_l^a(\phi')\phi'^*(\frac{\p}{\p u^a})
			+|Q\chi|^2 \psi'^a\otimes \phi'^*(\frac{\p}{\p u^a}) \\
		& +\frac{1}{3}\sum_{l,b}\big(\langle\psi'^b,\psi'^d\rangle\psi'^c
			-\langle\psi'^c,\psi'^b \rangle\psi'^d\big) \otimes
			\frac{\p \nu_l^b}{\p u^d}(\frac{\p \nu_l}{\p u^c})^{\top,a} \phi'^*(\frac{\p}{\p u^a})\\
		& -e_\al\cdot \na \phi'^a\cdot\chi^\al \otimes \phi'^*(\frac{\p}{\p u^a}).
	\end{split}
\end{equation}
In components, for each $a$,
\begin{equation}%
\label{eq-psi'-Componentwisely}
	\begin{split}
		\pd\psi'^a =\;&-\sum_{l,b}\na\phi'^d\cdot\psi'^b\frac{\p\nu_l^b}{\p u^d}\nu_l^a(\phi')+|Q\chi|^2 \psi'^a \\
		& +\frac{1}{3}\sum_{l,b}\big(\langle\psi'^b,\psi'^d\rangle\psi'^c
			-\langle\psi'^c,\psi'^b \rangle\psi'^d\big)
			\frac{\p \nu_l^b}{\p u^d}(\frac{\p \nu_l}{\p u^c})^{\top,a}
			-e_\al\cdot \na \phi'^a\cdot\chi^\al.
	\end{split}
\end{equation}
Here $\pd$ is the Dirac operator \(\pd\) on \(S\) and each $\psi'^a$ is a local pure spinor field.

Next we apply $\Hat{\phi}^*(f_*)$ to $\tau(\phi)$ to get
\begin{equation}
	\begin{split}
		\tau(\phi')-\sum_{\al}A(\phi)\big(T\phi (e_\al),T\phi (e_\al)\big)
		=\;& \frac{1}{2}\Hat{\phi}^*(f_*)\Rm^{\phi^*TN}(\psi,e_\al\cdot\psi)\phi_* e_\al
		-\frac{1}{12}\Hat{\phi}^*(f_*)(S\na R(\psi)) \\
		& -\Hat{\phi}^*(f_*) \big((\diverg V^j)\phi^*(\frac{\p}{\p y^j})+\na^{\phi^*TN}_{V^j} \phi^*(\frac{\p}{\p y^j})\big).
	\end{split}
\end{equation}
Since $\R^K$ is flat,
\begin{equation}
	LHS=\Delta\phi'-\sum_{\al} A(\phi')\big(T\phi' (e_\al),T\phi'(e_\al)\big)
	=\Delta \phi'
	+\sum_{\al}\frac{\p\phi'^a}{\p x^\al}\frac{\p\phi'^b}{\p x^\al}\frac{\p\nu_l^b}{\p u^a}(\phi') \phi^*(\nu_l).
\end{equation}
We deal with the terms on the right hand side as follows:
\begin{itemize}
	\item Using~\eqref{Rm-P} we get
		\begin{equation}
			\begin{split}
				\Hat{\phi}\big(&\Rm^{\phi^*TN}(\psi,e_\al\cdot\psi)\phi_* e_\al \big) \\
				&= \langle\psi^k,e_\al\cdot\psi^l\rangle \Rm(\frac{\p}{\p y^k}, \frac{\p}{\p y^l})T\phi (e_\al)\\
				&= \langle\psi^k,e_\alpha\left(\phi^j\right)\cdot\psi^l\rangle
				\Rm(\frac{\p}{\p y^k},\frac{\p}{\p y^l})\frac{\p}{\p y^j} \\
				&= \langle\psi^k,\na\phi^j\cdot\psi^l\rangle
				\big( P(A(\frac{\p}{\p y^j},\frac{\p}{\p y^l});\frac{\p}{\p y^k})
				-P(A(\frac{\p}{\p y^k},\frac{\p}{\p y^j});\frac{\p}{\p y^l})\big)  \\
				&= \langle\psi^k,\na\phi^j\cdot\psi^l\rangle P(A(\frac{\p}{\p y^j},\frac{\p}{\p y^l});\frac{\p}{\p y^k})+\langle\na\phi^j\cdot\psi^k,\psi^l\rangle P(A(\frac{\p}{\p y^k},\frac{\p}{\p y^j});\frac{\p}{\p y^l})\\
				&=2\langle\psi^k,\na\phi^j\cdot\psi^l\rangle  P(A(\frac{\p}{\p y^j},\frac{\p}{\p y^l});\frac{\p}{\p y^k})\\
				&=2\langle\psi^k,e_\al\cdot\psi^l\rangle  P\big(A(T\phi (e_\al),\frac{\p}{\p y^l});\frac{\p}{\p y^k}\big)
			\end{split}
		\end{equation}
		Hence
		\begin{equation}
			\begin{split}
				\frac{1}{2}\Hat{\phi}^*(f_*)\Rm^{\phi^*TN}(\psi,e_\al\cdot\psi)\phi_* e_\al
				=\;&\langle\psi^k,e_\al\cdot\psi^l\rangle \Hat{\phi}^*\left(f_*\right) P\big(A(T\phi(e_\al),\frac{\p}{\p y^l});\frac{\p}{\p y^k}\big)\\
				=\;&\langle\psi^k,e_\al\cdot\psi^l\rangle
				\frac{\p\phi'^a}{\p x^\al}\frac{\p u^b}{\p y^l}\frac{\p u^c}{\p y^k}
				\frac{\p\nu_l^b}{\p u^a}(\frac{\p\nu_l}{\p u^c})^\top  \\
				=\;&\langle\psi'^c,\na\phi'^a\cdot\psi'^b\rangle \frac{\p\nu_l^b}{\p u^a}(\frac{\p\nu_l}{\p u^c})^\top
			\end{split}
		\end{equation}
	\item To push $S\na R$ forward, we note that we can extend the local coordinate functions,
		which are defined in an open subset of $N$, so that they are constant in normal directions.
		Thus~$y^i$, $i=1,\dotsc,n$, are defined in a tubular neighborhood of a domain in $N$,
		which is an open subset of $\R^K$.
		The derivatives of $y^i$ with respect to $u^a$ are uniquely defined on $N$. Then
		\begin{equation}
			\begin{split}
				-\frac{1}{12}\Hat{\phi}^*(f_*)(S\na R(\psi))
				=\;&-\frac{1}{12}\Hat{\phi}^*(f_*)\big((\na\Rm)_{mjkl}\langle\psi^m,\psi^k\rangle\langle\psi^j,\psi^l\rangle\big) \\
				=\;&-\frac{1}{12}\big((\na\Rm)_{abcd}\langle\psi'^a,\psi'^c\rangle\langle\psi'^b,\psi'^d\rangle\big),
			\end{split}
		\end{equation}
		where
		\begin{equation}
			(\na\Rm)_{abcd}(x)=\big((\na\Rm)_{ijkl}\frac{\p y^i}{\p u^a}\frac{\p y^j}{\p u^b}
			\frac{\p y^k}{\p u^c}\frac{\p y^l}{\p u^d}\big)(\phi'(x)).
		\end{equation}
		Moreover, using Gauss equation again, one has
		\begin{equation}
			(\na \Rm)_{ijkl}=2(\la \na A_{ik}, A_{jl}\ra-\la \na A_{il}, A_{jk}\ra).
		\end{equation}
		where we have written $A_{ij}\equiv A(\frac{\p}{\p y^i}, \frac{\p}{\p y^j})$.
		See for example~\cite{branding2015some, jost2015geometric}. Hence,
		\begin{equation}
			\begin{split}
				-\frac{1}{12}\Hat{\phi}^*(f_*)&(S\na R(\psi)) \\
				&= -\frac{1}{6}\big(\la \na A_{ik}, A_{jl}\ra-\la \na A_{il}, A_{jk}\ra\big)
				\frac{\p y^i}{\p u^a}\frac{\p y^j}{\p u^b} \frac{\p y^k}{\p u^c}\frac{\p y^l}{\p u^d}
				\langle\psi'^a,\psi'^c\rangle\langle\psi'^b,\psi'^d\rangle   \\
				&\eqqcolon Z(A,\na A)_{abcd} \langle\psi'^a,\psi'^c\rangle\langle\psi'^b,\psi'^d\rangle.
			\end{split}
		\end{equation}

	\item In the same way as we have defined the vector fields $V^j,j=1,\dotsc,n$,
		we can define vector fields $V'^a, a=1\dotsc,K$, on $M$ by
		\begin{equation}
			\langle V'^a,W\rangle_{TM}=\langle e_\al\cdot W\cdot \chi^\al,\psi'^a\rangle_S,\quad\forall W\in\Gamma(TM).
		\end{equation}
		Then
		\begin{equation}
			\begin{split}
				\Hat{\phi}^*&(f_*) \big((\diverg V^j)\phi^*(\frac{\p}{\p y^j})+\na^{\phi^*TN}_{V^j} \phi^*(\frac{\p}{\p y^j})\big)\\
				&= e_\be\langle e_\al\cdot e_\be\cdot\chi^\al,\psi^j\rangle_S\frac{\p u^a}{\p y^j}\phi'^*(\frac{\p}{\p u^a})
				+V^j(\frac{\p u^a}{\p y^j}(\phi'))\phi'^*(\frac{\p}{\p u^a})
				-V^{j,\be}A(T\phi (e_\be),\frac{\p}{\p y^j}) \\
				&= e_\be(\langle e_\al\cdot e_\be\cdot\chi^\al,\psi'^a\rangle)\phi'^*(\frac{\p}{\p u^a})
				+\frac{\p \phi'^a}{\p x^\be} V'^{b,\be}\frac{\p \nu_l^b}{\p u^a}\phi'^*(\nu_l)\\
				&= (\diverg V'^a)\phi'^*(\frac{\p}{\p u^a})+\langle V'^b,\na\phi'^a\rangle\frac{\p\nu_l^b}{\p u^a}\phi'^*(\nu_l).
			\end{split}
		\end{equation}
\end{itemize}

Therefore the equation for $\phi'$ is
\begin{equation}%
\label{eq-phi'}
	\begin{split}
		\Delta \phi'
		=\;& -\sum_{\al,l}\frac{\p\phi'^a}{\p x^\al}\frac{\p\phi'^b}{\p x^\al}\frac{\p\nu_l^b}{\p u^a}(\phi') \phi'^*(\nu_l)
		+\sum_{b,l}\langle\psi'^c,\na\phi'^a\cdot\psi'^b\rangle
		\frac{\p\nu_l^b}{\p u^a}(\frac{\p\nu_l}{\p u^c})^\top(\phi')\\
		& +Z(A,\na A)_{abcd}\langle\psi'^a,\psi'^c\rangle\langle\psi'^b,\psi'^d\rangle
		-(\diverg V'^a)\phi'^*(\frac{\p}{\p u^a})-\langle V'^b,\na\phi'^a\rangle\frac{\p\nu_l^b}{\p u^a}\phi'^*(\nu_l).
	\end{split}
\end{equation}
In components, for each $a$,
\begin{equation}%
\label{eq-phi'-Componentwisely}
	\begin{split}
		\Delta \phi'^a
		=\;&-\sum_{\al,b,l}\frac{\p\phi'^c}{\p x^\al}\frac{\p\phi'^b}{\p x^\al}\frac{\p\nu_l^b}{\p u^c}(\phi')\nu_l^a(\phi')
		+\sum_{b,l}\langle\psi'^c,\na\phi'^d\cdot\psi'^b\rangle
		\frac{\p\nu_l^b}{\p u^d}(\frac{\p\nu_l}{\p u^c})^{\top,a} (\phi') \\
		& +Z^a(A,\na A)_{ebcd}\langle\psi'^e,\psi'^c\rangle\langle\psi'^b,\psi'^d\rangle
		-\diverg V'^a-\langle V'^b,\na\phi'^c\rangle\frac{\p\nu_l^b}{\p u^c}\nu_l^a(\phi').
	\end{split}
\end{equation}
As in~\cite{zhu2009regularity} and~\cite{chen2011boundary}, we shall transform the equation in a suitable form for later use.
Since $\phi'_* e_\al$ is tangent to $N$ while $\nu_l$ is perpendicular to $N$, they are orthogonal:
\begin{equation}\label{orthogonality}
	\sum_{b}\frac{\p \phi'^b}{\p x^\al} \nu_l^b=0, \quad \forall \al, \forall l.
\end{equation}
Hence
\begin{equation}
	\sum_{\al,b,l} \frac{\p\phi'^c}{\p x^\al} \frac{\p\phi'^b}{\p x^\al}\frac{\p\nu_l^a}{\p u^c} \nu_l^b=0,
\end{equation}
and we can add it to the first summand of~\eqref{eq-phi'-Componentwisely} to get a term of the form
\begin{equation}
	\sum_{\al,b,l} \frac{\p\phi'^b}{\p x^\al}\big(\frac{\p\phi'^c}{\p x^\al}\frac{\p\nu_l^a}{\p u^c} \nu_l^b
	-\frac{\p\phi'^c}{\p x^\al}\frac{\p\nu_l^b}{\p u^c}\nu_l^a \big)
	=\sum_{\al,b} \omega_\al^{ab} \frac{\p\phi'^b}{\p x^\al},
\end{equation}
with
\begin{equation}
	\omega_\al^{ab}=-\Big(\frac{\p\phi'^c}{\p x^\al}\frac{\p\nu_l^a}{\p u^c} \nu_l^b
	-\frac{\p\phi'^c}{\p x^\al}\frac{\p\nu_l^b}{\p u^c}\nu_l^a\Big)
	=-\omega_{\al}^{ba}.
\end{equation}
The second summand of~\eqref{eq-phi'-Componentwisely} can also be arranged into such a form.
Actually, using the symmetry~\eqref{symmetry of P}, we get
\begin{equation}
	\begin{split}
		&\sum_{b,l}\langle\psi'^c,\na\phi'^d\cdot\psi'^b\rangle\frac{\p\nu_l^b}{\p u^d}(\frac{\p\nu_l}{\p u^c})^{\top,a}\\
		=\;&
			\sum_{b,l}\langle\psi'^c,\na\phi'^b\cdot\psi'^d\rangle\frac{\p\nu_l^b}{\p u^d}(\frac{\p\nu_l}{\p u^c})^{\top,a}\\
		=\;&\frac{1}{2} \sum_{b,l}\big(\langle\psi'^c,\na\phi'^b\cdot\psi'^d\rangle
			+\langle\na\phi'^b\cdot\psi'^d,\psi'^c\rangle\big)
			\frac{\p\nu_l^b}{\p u^d}(\frac{\p\nu_l}{\p u^c})^{\top,a}\\
		=\;&\frac{1}{2} \sum_{\al,b,l} \langle\psi'^c,e_\al\cdot\psi'^d\rangle \frac{\p\phi'^b}{\p x^\al}
			\frac{\p\nu_l^b}{\p u^d}(\frac{\p\nu_l}{\p u^c})^{\top,a}
			+\langle e_\al\cdot\psi'^d,\psi'^c\rangle \frac{\p\phi'^b}{\p x^\al}
			\frac{\p\nu_l^b}{\p u^d}(\frac{\p\nu_l}{\p u^c})^{\top,a}\\
		=\;&\frac{1}{2}\sum_{\al,b,l}\langle\psi'^c,e_\al\cdot\psi'^d\rangle \frac{\p\phi'^b}{\p x^\al}
			\frac{\p\nu_l^b}{\p u^d}(\frac{\p\nu_l}{\p u^c})^{\top,a}
			-\langle\psi'^c,e_\al\cdot\psi'^d\rangle \frac{\p\phi'^b}{\p x^\al}
			\frac{\p\nu_l^b}{\p u^c}(\frac{\p\nu_l}{\p u^d})^{\top,a}.
	\end{split}
\end{equation}
Since $\phi'_* e_\al$ is tangent to $N$,
\begin{equation}
	\sum_b\frac{\p\phi'^b}{\p x^\al}\frac{\p\nu_l^b}{\p u^d}=\sum_b\frac{\p\phi'^b}{\p x^\al}(\frac{\p\nu_l}{\p u^d})^b
	=\sum_{b} \frac{\p\phi'^b}{\p x^\al}(\frac{\p\nu_l}{\p u^d})^{\top,b}.
\end{equation}
Thus the above term equals
\begin{equation}
	\begin{split}
		\frac{1}{2}\sum_{\al,b,l}\langle\psi'^c,e_\al\cdot\psi'^d\rangle
		\big((\frac{\p\nu_l}{\p u^d})^{\top,b}(\frac{\p\nu_l}{\p u^c})^{\top,a}
		-(\frac{\p\nu_l}{\p u^d})^{\top,a}(\frac{\p\nu_l}{\p u^c})^{\top,b}\big)
		\frac{\p\phi'^b}{\p x^\al}
		\equiv \sum_{\al,b} F^{ab}_\al \frac{\p\phi'^b}{\p x^\al},
	\end{split}
\end{equation}
with
\begin{equation}
	F^{ab}_\al=\sum_{\al,l} \langle\psi'^c,e_\al\cdot\psi'^d\rangle
	\big((\frac{\p\nu_l}{\p u^d})^{\top,b}(\frac{\p\nu_l}{\p u^c})^{\top,a}
	-(\frac{\p\nu_l}{\p u^d})^{\top,a}(\frac{\p\nu_l}{\p u^c})^{\top,b}\big)
	=-F^{ba}_\al.
\end{equation}

Similarly, using~\eqref{orthogonality}, the last summand of~\eqref{eq-phi'-Componentwisely} can be rearranged as
\begin{equation}
	\begin{split}
		\sum_{c,l} \langle V'^c,\na\phi'^b\rangle\frac{\p\nu_l^c}{\p u^b}\nu_l^a(\phi')
		&=\sum_{c,l,\al} \frac{\p \nu_l^c}{\p u^b}V'^c_\al\nu_l^a(\phi')\frac{\p\phi'^b}{\p x^\al} \\
		&=\sum_{b,c,l,\al}\big(\frac{\p \nu_l^c}{\p u^b}V'^c_\al\nu_l^a(\phi')-\frac{\p\nu_l^c}{\p u^a}V'^c_\al \nu^b_l(\phi')\big)
		\frac{\p\phi'^b}{\p x^\al} \\
		&\equiv -\sum_{\al,b}T^{ab}_\al\frac{\p\phi'^b}{\p x^\al}
	\end{split}
\end{equation}
where
\begin{equation}
	T^{ab}_\al=-\sum_{c}\left(\frac{\p \nu_l^c}{\p u^b}V'^c_\al\nu_l^a(\phi')-\frac{\p\nu_l^c}{\p u^a}V'^c_\al \nu^b_l(\phi')\right)
	=-T^{ba}_\al.
\end{equation}
\begin{rmk}
	Actually, for our proof of the local regularity of weak solutions, we don't need to write the second term
	and the last term into such an antisymmetric structure, see~\cite{branding2015some} for a similar treatment for a simpler model.
	Former regularity proofs, see e.g.~\cite{zhu2009regularity, wang2009regularity, chen2011boundary}, however, did need that structure.
	But it is also convenient to have such a structure.
\end{rmk}
Therefore, the equations for $\phi'$ appear in the elegant form:
\begin{equation}\label{new-eq-phi}
	\Delta\phi'^a=\sum_{b,\al}(\omega^{ab}_\al+F^{ab}_\al+T^{ab}_\al)\frac{\p\phi'^b}{\p x^\al}
	+Z^a(A,\na A)_{ebcd}\langle\psi'^e,\psi'^c\rangle\langle\psi'^b,\psi'^d\rangle
	-\diverg V'^a,
\end{equation}
for $a=1,\dotsc,K$, where the coefficients of first derivative of $\phi'$ are antisymmetric.

\section{Regularity of weak solutions}

We now come to the crucial contribution of our paper, the regularity of weak solutions.
In order to make the action functional $\A$ well-defined and finite-valued, we need to assume
\begin{align}
	\phi &\in W^{1,2}(M,N), &
	\psi &\in W^{1,4/3}(\Gamma(S\otimes \phi^*TN)).
\end{align}
The issue then is higher regularity of such weak solutions.
More precisely, we shall show that~$(\phi,\psi)$ are smooth when they satisfy~\eqref{EL-eq} in the weak sense.
By the Sobolev embedding theorem, $\phi \in L^p(M,N)$ for any $p\in [1,\infty)$ and $\psi\in L^4(\Gamma(S\otimes\phi^*TN))$.
Since $f\colon N\to \R^K$ is a smooth embedding, $(\phi',\psi')$ have the same regularity as $(\phi,\psi)$,
and so it suffices to show smoothness of the former.

As the regularity is a local issue, we can take $\phi'\colon B_1\to\mathbb{R}^K$ defined in the euclidean unit disc $B_1\subset \R^2\cong \mathbb{C}^1$.
Over \(B_1\) the bundle \(S\otimes \phi^*TN'\) is trivial with typical fiber \(\mathbb{C}^2\otimes\mathbb{R}^K\).
Hence \(\psi'\colon B_1\to \mathbb{C}^2\otimes\mathbb{R}^K\) is a vector valued function.

\subsection{}
As we have seen, $\psi'$ satisfies~\eqref{eq-psi'} or equivalently~\eqref{eq-psi'-Componentwisely}.
By the following lemma, which will be proved in Section~\ref{Sec:ProofOfLemma}, all powers of $\psi'$ are integrable.
\begin{lemma}\label{reg for Dirac}
	Let $p\in(4,\infty)$ and $\vph\in L^4(B_1,\mathbb{C}^2\otimes \R^K)$ be a weak solution of the nonlinear system
	\begin{equation}
		\pd \vph^i =A^i_j \vph^j +B^i, \quad 1\le i \le K,
	\end{equation}
	where $A\in L^2(B_1, \mathfrak{gl}(2, \mathbb{C})\otimes\mathfrak{gl}(K, \mathbb{R}))$ and $B\in L^2(B_1,\mathbb{C}^2\otimes\R^K)$.
	There exists a $\varepsilon_0=\varepsilon_0(p)>0$ such that if $\|A\|_{L^2(B_1)}\le \varepsilon_0$, then $\vph\in L^p_{loc}(B_1)$.
\end{lemma}

It follows from Lemma~\ref{reg for Dirac} that $\psi'\in L^p_{loc}(B_1)$ for any $p\in[1,\infty)$.
Since locally the Dirac operator is given by the classical Cauchy--Riemann operators \(\partial_{z}\) and \(\partial_{\overline{z}}\), it follows from the elliptic theory that $\psi'\in W^{1,q}(B_{1/2})$ for any $q\in [1,2)$.

\subsection{}
We use the aforementioned Rivi\`ere's regularity theory to deal with $\phi'$.
More precisely, we use the following result which is an extension of~\cite{riviere} to improve the regularity of $\phi'$.

\begin{thm}\label{Sharp-Topping}
	\emph{(\cite{riviere2010conformally, sharp2013decay})}
	Let $p\in (1,2)$. Suppose that $u\in W^{1,2}(B_1,\R^K)$ is a weak solution of
	\begin{equation}
		-\Delta u=\Omega \na u+f,
	\end{equation}
	where $\Omega \in L^2(B_1,\mathfrak{so}(K)\otimes \R^2)$ and $f\in L^p(B_1,\R^K)$. Then $u\in W^{2,p}_{loc}(B_1)$.
\end{thm}

In the previous section we have written the equation for $\phi'$ into such a form, see~\eqref{new-eq-phi}.
Since we have seen
$\psi'\in L^p_{loc}(B_1), 1\le p<\infty$, the hypotheses of Theorem~\ref{Sharp-Topping} are satisfied.
Thus we can conclude that $\phi'\in W^{2,p}_{loc}(B_1)$ for any $p\in[1,2)$.
It follows from the Sobolev embedding theorems that $\phi'\in W^{1,q}(B_{1/2})$ for any $q\in [1,\infty)$.

\subsection{}
We can now apply the standard elliptic theory for a bootstrap argument,
see e.g.~\cite{begehr1994complex,gilbarg2001elliptic}, and hence conclude that
$(\phi',\psi')$ are smooth. The smoothness of $\phi$ then follows directly.
For $\psi$, one can use~\eqref{local form of EL-eq on N} and the elliptic theory for Cauchy--Riemann operators (e.g.~\cite{begehr1994complex}) to conclude that $\psi$ is also smooth.
Therefore the full regularity of weak solutions is obtained, completing the proof of Theorem~\ref{theorem 1}.

\section{Proof of Lemma~\ref{reg for Dirac}}%
\label{Sec:ProofOfLemma}

In this section, we provide the proof of Lemma~\ref{reg for Dirac}.
We shall use the Dirac type equation to improve the integrability of the spinor.
Results of this type were first obtained by~\cite{wang2010remark} and further developed in~\cite{sharp2016regularity, branding2015some}.
Actually a stronger result holds in general.
Before stating the general result, we recall some basic facts on Morrey spaces, see for example~\cite{giaquinta1983multiple}.

Let $U$ be a domain in $\R^n$. For $0\le \lambda\le n$ and $1\le p<\infty$, the Morrey space on $U$ is defined as
\begin{equation}
	\MS{p,\lambda}(U)\coloneqq\Big\{u\in L^p(U)\big| \|u\|_{\MS{p,\lambda}(U)}<\infty \Big\}.
\end{equation}
Here the $(p,\lambda)$-Morrey norm of $u$ is defined by
\begin{equation}
	\|u\|_{\MS{p,\lambda}(U)}\coloneqq\sup_{x\in U,\;r>0} \Big(\frac{r^\lambda}{r^n}\int_{B_r(x)\cap U} |u(y)|^p \dd y \Big)^{1/p}.
\end{equation}
Note that on a bounded domain $U\subset \R^n$, for $1\le p <\infty$ and $0\le \lambda\le n$, it holds that
\begin{equation}
	L^{\infty}(U)=\MS{p,0}(U)\subset \MS{p,\lambda}(U) \subset \MS{p,n}(U)= L^p(U).
\end{equation}
In this section we consider a map $\vph\colon B_1\to \mathbb{C}^L\otimes \R^K$ satisfying a first order elliptic system, where $B_1\subset \R^n$ is the euclidean unit ball and $\mathbb{C}^L\otimes \R^K$ is supposed to be the typical fiber of a twisted complex spinor bundle over \(B_1\).

\begin{lemma}\label{Morrey type regularity}
	Let $n\ge 2$ and $4<p<+\infty$.
	Let $\vph\in \MS{4,2}(B_1, \mathbb{C}^L\otimes \R^K)$ be a weak solution of the nonlinear system
	\begin{equation}\label{nonlinear dirac system-2}
		\pd \vph^i=A^i_j \vph^j + B^i, \quad 1\le i\le K,
	\end{equation}
	where $A\in \MS{2,2}(B_1, \mathfrak{gl}(L,\mathbb{C})\otimes\mathfrak{gl}(K,\R))$ and $B \in \MS{2,2}(B_1,\mathbb{C}^L\otimes \R^K)$.
	There exists $\varepsilon_0=\varepsilon_0(n,p)>0$ such that if
	\begin{equation}
		\|A\|_{\MS{2,2}(B_1)}\le \varepsilon_0,
	\end{equation}
	then $\vph\in L^p_{loc}(B_1)$. Moreover, for any $U\Subset B_1$,
	\begin{equation}
		\|\vph\|_{L^p(U)} \le C(n,p,U)\big(\|\vph\|_{\MS{4,2}(B_1)}+\|B\|_{\MS{2,2}(B_1)}\big).
	\end{equation}
\end{lemma}
The proof is motivated from that in~\cite{wang2010remark} and is adapted to this system with minor changes.
The idea is to use the fundamental solution of the Euclidean Dirac operator and apply Riesz potential estimates.
Thanks to the Bochner-Lichnerowicz-Weitzenb\"ock type formulas, e.g.\ see~\cite[Theorem II.8.17]{lawson1989spin},~\cite[Lemma 4.1]{tolksdorf2001clifford},~\cite[Theorem 4.4.2]{jost2008riemannian}, the fundamental solution of the Euclidean Dirac operator can be derived from that of the Euclidean Laplacian.
We remark that the $\MS{2,2}$-assumption on $B$ here fits quite well to the proof.
\begin{proof}
	Applying $\pd$ to~\eqref{nonlinear dirac system-2}, we have, for $1\le i\le K$,
	\begin{equation}
		-\Delta \vph^i=\pd^2 \vph^i=\pd(A^i_j\vph^j+B^i)
	\end{equation}
	in the sense of distributions.

	Let $x_0\in B_1 $, $|x_0|<1$, and let $0<R<1-|x_0|$.
	Take a cutoff function $\eta\in C_0^\infty(B_{R}(x_0))$ such that $0 \le \eta \le 1$ and $\eta\equiv 1$ on $B_{ R/2}(x_0)$.
	For each $1\le i\le K$, define $g^i\colon \R^n \to \mathbb{C}^L$ by
	\begin{equation}\label{good part-2}
		g^i(x) =\int_{\R^n} \frac{\p G(x,y)}{\p y^\al} \frac{\p}{\p y^\al} \cdot
		\big(\eta^2( A^i_j\vph^j+B^i)\big)(y) \dd y
	\end{equation}
	where $G(x,y)$ is the fundamental solution of $\Delta$ on $\R^n$.
	Thus
	\begin{equation}
		\begin{split}
			-\Delta g^i &= \pd \big(\eta^2 (A^i_j\vph^j+B^i)\big) \\
			&= \pd \big( A^i_j\vph^j+B^i\big) \quad \quad \text{in \(B_{R/2}(x_0)\)}.
		\end{split}
	\end{equation}
	Setting $h^i\coloneqq\vph^i-g^i$, we see that $h^i$, $1\le i\le K$, are harmonic in $B_{R/2}(x_0)$:
	\begin{equation}
		\Delta h^i=0 \quad \text{in \(B_{R/2}(x_0)\)}.
	\end{equation}
	Note that
	\begin{equation}
		|g^i(x)|\le C\int_{\R^n}\frac{1}{|x-y|^{n-1}}(\eta^2 |A^i_j\vph^j+B^i|) \dd y=C I_1(\eta^2 (A\vph+B)),
	\end{equation}
	where $I_1$ is the Riesz potential operator. By Adams' inequality~\cite[Theorem 3.1]{adams1975note}, for $1<q<\lambda\le n$,
	\begin{equation}
		\|I_1(\eta^2 (A\vph+B))\|_{\MS{\frac{\lambda q}{\lambda-q},\lambda}(\R^n)}
		\le C\|\eta^2 |A\vph+B|\|_{\MS{q,\lambda}(\R^n)}.
	\end{equation}

	\textbf{\underline{Step 1:}} By hypothesis we have
	\begin{equation}
		\begin{split}
			\|\eta^2(A\vph+B)\|_{\MS{\frac{4}{3},2}(\R^n)}
			& \le \|(\eta A)(\eta\vph)\|_{\MS{\frac{4}{3},2}(\R^n)}+\|\eta^2 B\|_{\MS{\frac{4}{3},2}(\R^n)} \\
			& \le \|\eta A\|_{\MS{2,2}(\R^n)} \|\eta \vph\|_{\MS{4,2}(\R^n)} +\|\eta^2 B\|_{\MS{\frac{4}{3},2}(\R^n)} \\
			& \le \|A\|_{\MS{2,2}(B_R(x_0))} \|\vph\|_{\MS{4,2}(B_R(x_0))}+\|B\|_{\MS{\frac{4}{3},2}(B_R(x_0))}\\
			& \le \|A\|_{\MS{2,2}(B_R(x_0))} \|\vph\|_{\MS{4,2}(B_R(x_0))}+C R^{\frac{1}{2}}\|B\|_{\MS{2,2}(B_R(x_0))}.
		\end{split}
	\end{equation}
	With $q=\frac{4}{3}, \lambda=2, \frac{\lambda q}{\lambda-q}=4$. We get
	\begin{equation}
		\begin{split}
			\|g\|_{\MS{4,2}(\R^n)}&\le C\|I_1(\eta^2 (A\vph+B))\|_{\MS{4,2}(\R^n)}
			\le C\|\eta^2(A\vph+B)\|_{\MS{\frac{4}{3},2}(\R^n)} \\
			&\le C\varepsilon_0 \|\vph\|_{\MS{4,2}(B_R(x_0))}+C R^{\frac{1}{2}}|B|.
		\end{split}
	\end{equation}
	where we have denoted $|B|\equiv \|B\|_{\MS{2,2}(B_1)}$.

	Note that $|h^i|^4$ is subharmonic in $B_{R/2}(x_0)$:
	\begin{equation}
		\Delta |h^i|^4=\Delta{(h^i\overline{h^i})}^2=2|\na(h^i\overline{h^i})|^2+2|h^i|^2\big((\Delta h^i)\overline{h^i}+2|\na h^i|^2+ h^i\overline{\Delta h^i}\big)\ge0
	\end{equation}
	since $\Delta h^i=0$.
	Hence $\fint_{B_r(x)}|h^i|^4\dd y$ is a nondecreasing function in $r$, which implies that for any $1\le i\le m$ and any $\theta\in (0,1/6)$,
	\begin{equation}
		\|h^i\|_{\MS{4,2}(B_{\theta R}(x_0))} \le {(4\theta)}^{1/2}\|h^i\|_{\MS{4,2}(B_{R/2}(x_0))}.
	\end{equation}
	Recalling $\vph^i=g^i+h^i$, we get
	\begin{equation}
		\begin{split}
			\|\vph\|_{\MS{4,2}(B_{\theta R}(x_0))}
			&\le \|g\|_{\MS{4,2}(B_{\theta R}(x_0))}+\|h\|_{\MS{4,2}(B_{\theta R}(x_0))} \\
			&\le  C\varepsilon_0 \|\vph\|_{\MS{4,2}(B_{R}(x_0))} +C|B| R^{\frac{1}{2}}
			+2\theta^{1/2}\|h\|_{\MS{4,2}(B_{R/2}(x_0))} \\
			&\le C\varepsilon_0 \|\vph\|_{\MS{4,2}(B_{R}(x_0))} +C|B| R^{\frac{1}{2}}
			+2\theta^{1/2}\big(\|\vph\|_{\MS{4,2}(B_{R/2}(x_0))}+\|g\|_{\MS{4,2}(B_{R/2}(x_0))} \big) \\
			&\le C_0(\varepsilon_0+\theta^{1/2})\|\vph\|_{\MS{4,2}(B_{R}(x_0))}+C|B|R^{\frac{1}{2}}.
		\end{split}
	\end{equation}
	Fix any $\be\in (0,\frac{1}{2})$, we can find a $\theta\in(0,\frac{1}{2})$ such that $2C_0 \theta^{1/2} \le \theta^\be$.
	Then take $\varepsilon_0$ small enough such that $2C_0 \varepsilon_0\le \theta^\be$. With such a choice we have
	\begin{equation}\label{iteration inequality-2}
		\|\vph\|_{\MS{4,2}(B_{\theta R}(x_0))}\le \theta^\be \|\vph\|_{\MS{4,2}(B_{R}(x_0))}+C|B|R^{\frac{1}{2}}.
	\end{equation}
	Note that~\eqref{iteration inequality-2} holds for any $0<R<1-|x_0|$.
	Thus we can start the following iteration procedure.

	Let $R<1-|x_0|$. Then for any $0<r< R$, there exists a unique $k\in\mathbb{N}$ such that
	$\theta^{k+1}R <r\le \theta^k R$. (The case $k=0$ is trivial, and we may thus assume $k\ge 1$). Hence we have
	\begin{equation}
		\begin{split}
			\|\vph\|_{\MS{4,2}(B_{r}(x_0))}
			& \le \|\vph\|_{\MS{4,2}(B_{\theta^k R}(x_0))}
			\le \theta^\be\|\vph\|_{\MS{4,2}(B_{\theta^{k-1}R}(x_0))}+C|B|(\theta^{k-1}R)^{\frac{1}{2}} \\
			& \le \theta^{2\be} \|\vph\|_{\MS{4,2}(B_{\theta^{k-2}R}(x_0))}
			+C|B|[\theta^\be(\theta^{k-2}R)^{\frac{1}{2}}+(\theta^{k-1}R)^{\frac{1}{2}}]\\
			& \le \theta^{k\be} \|\vph\|_{\MS{4,2}(B_{R}(x_0))}
			+C|B|R^{\frac{1}{2}}\theta^{(k-1)\be}[1+\theta^{\frac{1}{2}-\be}+\cdots+\theta^{(\frac{1}{2}-\be)(k-1)}]\\
			& \le \frac{1}{\theta^\be}  \theta^{(k+1)\be} \|\vph\|_{\MS{4,2}(B_{R}(x_0))}
			+\frac{C|B|R^{\frac{1}{2}-\be}}{\theta^{2\be}}
			\frac{1-\theta^{(\frac{1}{2}-\be)k}}{1-\theta^{\frac{1}{2}-\be}} (\theta^{k+1}R)^\be \\
			& \le \frac{1}{\theta^\be}  \big(\frac{r}{R}\big)^\be \|\vph\|_{\MS{4,2}(B_{R}(x_0))}
			+\frac{C|B|}{\theta^{2\be}-\theta^{\frac{1}{2}+\be}} r^\be
		\end{split}
	\end{equation}
	where we used $R\le 1$ in the last inequality.
	In particular this implies that
	\begin{equation}
		\Big(\frac{1}{r^{n-2+4\be}}\int_{B_r(x_0)}|\vph|^4\dd y\Big)^{\frac{1}{4}}
		\le \frac{1}{(\theta R)^{\be}}  \|\vph\|_{\MS{4,2}(B_1)} +\frac{C|B|}{\theta^{2\be}-\theta^{\frac{1}{2}+\be}}
	\end{equation}
	If we restrict to $|x_0|<\frac{1}{4}$ and $R=\frac{1}{2}$, we see that $\vph \in \MS{4,2-4\be}(B_{\frac{1}{4}})$, with
	\begin{equation}
		\|\vph\|_{\MS{4,2-4\be}(B_{1/4})} \le C \|\vph\|_{\MS{4,2}(B_1)}+C\|B\|_{\MS{2,2}(B_1)}.
	\end{equation}
	for some universal constant $C=C(n,\be)$.

	\textbf{\underline{Step 2:}} We improve the integrability. Let $|x_0|<\frac{1}{4}$ and $0<R<\frac{1}{4}-|x_0|$.
	Take a cutoff function $\eta \in C_0^\infty(B_R(x_0))$ and define $g^i, h^i$ as before.
	Note that
	\begin{equation}
		\begin{split}
			\|\eta^2 (A\vph+B)\|_{\MS{\frac{4}{3}, 2-\frac{4\be}{3}}(\R^n)}
			& \le \|\eta^2 A\vph\|_{\MS{\frac{4}{3}, 2-\frac{4\be}{3}}(\R^n)}
			+\|\eta^2 B\|_{\MS{\frac{4}{3}, 2-\frac{4\be}{3}}(\R^n)} \\
			& \le \|\eta A\|_{\MS{2,2}(\R^n)} \|\eta\vph\|_{\MS{4,2-4\be}(\R^n)}
			+\|\eta^2 B\|_{\MS{\frac{4}{3}, 2-\frac{4\be}{3}}(\R^n)} \\
			& \le \|A\|_{\MS{2,2}(B_R(x_0))} \|\vph\|_{\MS{4,2-4\be}(B_R(x_0))}
			+\| B\|_{\MS{\frac{4}{3}, 2-\frac{4\be}{3}}(B_R(x_0))} \\
			& \le \|A\|_{\MS{2,2}(B_R(x_0))} \|\vph\|_{\MS{4,2-4\be}(B_R(x_0))}+C\| B\|_{\MS{2,2}(B_R(x_0))}R^{\frac{1}{2}-\be}\\
			& \le \varepsilon_0\|\vph\|_{\MS{4,2-4\be}(B_R(x_0))}+C|B|R^{\frac{1}{2}-\be}.
		\end{split}
	\end{equation}
	With $q=\frac{4}{3}$ and
	$\lambda=2-\frac{4\be}{3}$, (note that we need $1<q<\lambda\le n$, which requires $\be<\frac{1}{2}$), we see that,
	$\frac{\lambda q}{\lambda-q}=\frac{4(3-2\be)}{3-6\be}$, and
	\begin{equation}
		\begin{split}
			\|g^i\|_{\MS{\frac{4(3-2\be)}{3-6\be},2-\frac{4\be}{3}}(\R^n)}
			&  \le C\|I_1(\eta^2 (A\vph+B))\|_{\MS{\frac{4(3-2\be)}{3-6\be},2-\frac{4\be}{3}}(\R^n)}
			\le C\|\eta^2 (A\vph+B)\|_{\MS{\frac{4}{3}, 2-\frac{4\be}{3}}(\R^n)} \\
			&  \le C\varepsilon_0\|\vph\|_{\MS{4,2-4\be}(B_R(x_0))}+C|B|R^{\frac{1}{2}-\be}.
		\end{split}
	\end{equation}
	Again, $h^i$ is harmonic in $B_{R/2}(x_0)$ in the sense of distributions and $h^i\in L^4(B_{R/2}(x_0))$,
	by Weyl's lemma, it is smooth in $B_{R/2}(x_0)$, see e.g.~\cite[Corollary 1.2.1]{jost2013partial}.
	By shrinking the radius $R$ a little, we may assume $h^i\in L^{\infty}(B_{R/2}(x_0))$.
	Actually, by Harnack inequality in the disk together with mean value equality one has, for any $R'<R$,
	\begin{equation}
		\begin{split}
			\|h\|_{L^{\infty}(B_{R'/2}(x_0))}
			\le & C(R,R',n) |h(x_0)| \le C(R,R',n)\|h\|_{L^1(B_{R'/2}(x_0))} \\
			\le & C\big(\|\vph\|_{\MS{4,2}(B_1)} +\|B\|_{\MS{2,2}(B_1)} \big).
		\end{split}
	\end{equation}
	Thus if we restrict to $|x_0|\le \frac{1}{16}$ and $R=\frac{1}{8}$, we see that
	$h^i\in \MS{\frac{4(3-2\be)}{3-6\be},2-\frac{4\be}{3}}(B_{\frac{1}{16}})$.
	By elliptic theory,
	$\|h\|_{\MS{\frac{4(3-2\be)}{3-6\be},2-\frac{4\be}{3}}(B_{\frac{1}{16}})}$
	can be controlled by $\|\vph\|_{\MS{4,2}(B_1)}$ and $\|B\|_{\MS{2,2}(B_1)}$.

	Finally recall that
	\begin{equation}
		\vph^i=h^i+g^i.
	\end{equation}
	It follows that
	\begin{equation}\label{higher Morrey norm of psi-2}
		\begin{split}
			\|\vph\|_{\MS{\frac{4(3-2\be)}{3-6\be},2-\frac{4\be}{3}}(B_{\frac{1}{16}})}
			\le &  \|g\|_{\MS{\frac{4(3-2\be)}{3-6\be},2-\frac{4\be}{3}}(B_{\frac{1}{16}})}
			+ \|h\|_{\MS{\frac{4(3-2\be)}{3-6\be},2-\frac{4\be}{3}}(B_{\frac{1}{16}})} \\
			\le & C(\be,n)\big(\|\vph\|_{\MS{4,2}(B_1)}+\|B\|_{\MS{2,2}(B_1)}\big).
		\end{split}
	\end{equation}

	\textbf{\underline{Step 3:}} We note that~\eqref{higher Morrey norm of psi-2} holds for any give $0<\be<\frac{1}{2}$.
	Since
	\begin{equation}
		\lim_{\be \nearrow \frac{1}{2}} \frac{4(3-2\be)}{3-6\be}=+\infty,
	\end{equation}
	and
	\begin{equation}
		\MS{\frac{4(3-2\be)}{3-6\be},2-\frac{4\be}{3}}(B_{\frac{1}{16}})
		\hookrightarrow L^{\frac{4(3-2\be)}{3-6\be}}(B_{\frac{1}{16}}),
	\end{equation}
	we conclude that $\vph\in L^p(B_{\frac{1}{16}})$ for any $4<p<+\infty$ and
	\begin{equation}
		\|\vph\|_{L^p(B_{\frac{1}{16}})} \le C(n,p)\big(\|\vph\|_{\MS{4,2}(B_1)}+\|B\|_{\MS{2,2}(B_1)}\big).
	\end{equation}
	This completes the proof of the Lemma.
\end{proof}

Finally note that in the 2-dimensional case,
\begin{equation}
	\MS{2,2}(B_1)=L^2(B_1), \quad \MS{4,2}(B_1)=L^4(B_1).
\end{equation}
So Lemma~\ref{reg for Dirac} follows from Lemma~\ref{Morrey type regularity}.

\end{document}